\documentclass[journal,twoside,web]{IEEEtran}
\usepackage{amsthm} 

\usepackage{cite}
\usepackage{amsmath,amssymb,amsfonts}
\usepackage{algorithmic}
\usepackage{graphicx}
\usepackage{algorithm,algorithmic}
\usepackage{hyperref}

\usepackage{mathrsfs}
\usepackage{siunitx}
\usepackage{pgfplots}
\usepackage{tikz}
 \pgfplotsset{compat=newest}
  \usetikzlibrary{plotmarks}
  \usetikzlibrary{arrows.meta}
  \usetikzlibrary{shapes.geometric,calc}
  \usepgfplotslibrary{patchplots}

\newcommand{\legendline}[1]{%
  \begin{tikzpicture}[baseline=-0.6ex]
    \draw[#1, line width=1.2pt] (0,0) -- (0.5,0);
  \end{tikzpicture}%
}
\newcommand{\legendstar}[1]{%
  \begin{tikzpicture}[baseline=-0.6ex]
    \node[star,star points=5,star point ratio=2.25,
          draw=#1,fill=#1,minimum size=5pt,inner sep=0pt]{};
  \end{tikzpicture}%
}
  
  \usepackage{grffile}

\usepackage{subcaption}

\newlength\figH
\newlength\figW
\usepackage{url}

\makeatletter
\let\NAT@parse\undefined
\makeatother
\allowdisplaybreaks

\let\ieeebibliography\thebibliography

\usepackage[numbers,sort&compress]{natbib}

\renewcommand\thebibliography[1]{\ieeebibliography{#1}}

\usepackage{hyperref}

\usepackage[commentmarkup=uwave,authormarkup=none]{changes}
\definechangesauthor[color=purple]{TC} 
\definechangesauthor[color=blue]{AC} 
\definechangesauthor[color=olive]{JO}
\definechangesauthor[color=blue]{rev1}
\definechangesauthor[color=red]{rev1com}

\usepackage[IEEE]{altthm-styles}

\DeclareMathOperator{\interior}{int}

\newcommand{\CS}{\mathcal{X}}       
\newcommand{\TS}{\mathcal{X}_T}       
\newcommand{\RA}[2]{R_{[#1]}(#2)} 
\newcommand{\Ufcnset}[1]{\mathscr U_{[#1]}} 
\newcommand{\uu}{\mathbf{u}}
\newcommand{\xx}{\mathbf{x}}
\newcommand{\Uset}{\mathcal{U}}
\newcommand{\Tbdl}[1]{\mathcal{T}_{#1}}

\newcommand{\ip}[1]{\langle#1\rangle}
\newcommand{\dd}{\mbox{d}}

\newcommand{\Feas}{\mathcal{R}}
\newcommand{\ie}{i.e.\;}

\newcommand{\clf}[1][\empty]{\ifx#1\empty V \else V(#1)\fi}
\newcommand{\cbf}[1][\empty]{\ifx#1\empty h \else h(#1)\fi}

\newcommand{\Vap}{\hat V}
\newcommand{\Rap}{\hat{\mathcal{R}}}
\newcommand{\hap}{\hat h}

\newcommand{\uuinf}{\kappa_\partial}
\newcommand{\uuinfhat}{\hat{\kappa}_\partial}

\usepackage{textcomp}
\def\BibTeX{{\rm B\kern-.05em{\sc i\kern-.025em b}\kern-.08em
    T\kern-.1667em\lower.7ex\hbox{E}\kern-.125emX}}
\begin{document}

\title{Safe-by-Design: Approximate Nonlinear Model Predictive Control with Real Time Feasibility$^\star$}

\author{Jan Olucak, Arthur Castello B. de Oliveira, and Torbjørn Cunis, \IEEEmembership{Member,~IEEE} 
\thanks{$^\star$This research is partially supported by the Ministry of Science, Research and Arts of the state of Baden-Württemberg under funding number MWK32-7531-49/13/7 for the project DaSO: Data-driven Spacecraft Operations.}
\thanks{JO and TC are with the University of Stuttgart, 70569 Stuttgart, Germany~\\
        {\tt\small \{ jan.olucak | tcunis \} @ifr.uni-stuttgart.de}}%
\thanks{ACB is with Northeastern University, Boston, MA, USA~\\
        {\tt\small a.castello@northeastern.edu}}%
}

\maketitle

\begin{abstract}
 This paper proposes a computationally lightweight, continuous-time receding-horizon nonlinear model predictive control (MPC) approach with infinitesimal prediction horizon.  Unlike discrete-time {one-step} MPC schemes, this methodology permits the formulation of  small-sized convex quadratic programs (QPs)
for feedback {yet recovers the} theoretical guarantees  {of} quasi-infinite horizon MPC. The convex QP comes with a small computational footprint, which is advantageous for real-{time} application where runtime guarantees must be given. We demonstrate the effectiveness of the proposed approach when compared to other constrained control techniques through numerical experiments for nonlinear constrained spacecraft control.


\end{abstract}


\section{Introduction}
\label{sec:introduction}
\IEEEPARstart{M}{odel predictive control} (MPC)~\cite{gruenePannek2017} is a well-established control approach for a wide range of dynamical systems under constraints \cite{eren_model_2017,yu2021,nguyen2021}, providing a unified framework that ensures both safety and asymptotic stability. Based on an optimal control problem, MPC accommodates nonlinear system dynamics whilst handling input, state, and (complex) path constraints. However, if the computational demands cannot be met for complex, nonlinear or large-scale problems, the safety of MPC is compromised.  {The major goal of this work is to establish a real-time feasible  {nonlinear} MPC law  {which both provides} a low computational footprint and  {preserves the} strong theoretical  {guarantees of MPC}.}

To address the issue of computation time in MPC, various strategies have been discussed:
Explicit MPC~\cite{alessio2009,grancharova2012} precomputes the optimal control law offline, which can then be evaluated online in real-time using lookup tables. However, exponential growth of the precomputation with the number of system states makes this approach impractical for realistic applications. 
Alternatively, time-distributed optimization \cite{liao-mcpherson_time-distributed_2020} limits the number of iterations to solve the optimal control problem at each sampling time. This approach enables real-time implementation but compromises optimality and feasibility. Consequently, traditional MPC stability guarantees may no longer apply \cite{leung2021}. 
Other suboptimal MPC methods include the real-time iteration (RTI) scheme~\cite{diehl_real-time_2005}, advanced step MPC~\cite{zavala2009}, or combinations of both \cite{nurkanovic2019,nurkanovic2020}. 
To guarantee feasibility and asymptotic stability with low computational effort, some approaches rely on  precomputed sequences of invariant sets to reduce the prediction horizon to a single step for linear and nonlinear systems~\cite{angeli_ellipsoidal_2002,limon_robust_2003}. Here, { {asymptotic stability} is enforced by an additional constraint that requires an auxiliary optimization to be solved before evaluating the MPC feedback law. Our previous approach~\cite{olucak_nonlinear_2024} follows a similar philosophy but avoids the additional contraction constraint and instead augments the terminal penalty to achieve closed-loop asymptotic stability.  {More recently, the concept of nonlinear tracking MPC with artificial references was introduced in~\cite{limon2018a}. In case of  {a} reference-dependent terminal penalty and set, the prediction horizon can be reduced to arbitrarily small horizons~\cite{kohler2020c,kohler2020b,kohler2024},  {possibly achieving} a prediction horizon equal to one  {as well}. Recursive feasibility is ensured by the artificial reference and compared to classical MPC with terminal conditions, the reference-dependent terminal ingredients are not restricted to {a single equilibrium }point. However, even for a horizon equal to one, the underlying optimal control problem is still nonlinear and additional decision variables and constraints  {lead to increased complexity.} Furthermore, the set of reachable setpoints (output space) must be convex, which prohibits its direct application in motion planning problems with collision avoidance and needs special treatment~\cite[Section~IV.A]{krupa2024}.}
The computational effort of MPC can further be reduced by warm-starting with a previous solution~\cite{gruenePannek2017}. Additionally, using more dedicated warm-starting techniques~\cite{zeilinger_real-time_2014} or efficient nonlinear programming solvers, such as those in~\cite{zanelli2020,pas2022a,vanroye2023fatrop}, can be advantageous.

 {Except for explicit MPC, all of the above approaches result in nonlinear{ and usually} nonconvex parameter optimization problems, for which no runtime guarantee can be given.  {E}ven with specialized real-time schemes (e.g., time-distributed approaches) the computational effort might be still too high for applications with severe {ly limited} computational resources  {such as} space applications.
On the contrary, convex optimization problems are easier to solve and theoretical upper bounds on the number of iterations can be given~\cite[Subsection~11.5.3]{boyd2004}  {under mild assumptions}, which is advantageous for real- {time} applications.}

Like MPC techniques, strategies involving control Lyapunov \cite{isidori_nonlinear_1995} and control barrier \cite{wieland2007} functions (CLF/CBF) are also very popular for solving constrained nonlinear control problems.  {In particular, CLF and CBF constraints allow for efficient convex quadratic program (QP) formulations~\cite{li_survey_2023,garg2024}, and enable online computation of safe and stabilizing control laws. For example, CLFs can be used to find a stabilizing controller which minimizes control effort~\cite{freeman_inverse_1996}, whereas CBFs can be used to ensure safety through a \textit{safety filter}~\cite{gurriet_towards_2018}, or CLF and CBF(s) can be  {combined} \cite{ames_control_2017, schneeberger_advanced_2026}. Despite their excellent computational properties,  {{QP}-based} CLF/CBF approaches cannot account for state costs and stability must be ensured through constraints. The latter {is sometimes relaxed} to ensure feasibility,  {where} a slack variable is added to the QP to trade off between safety (CBF) and stability (CLF).
CLF and CBFs have also been incorporated into  discrete-time MPC frameworks
to establish a horizon-one MPC ~\cite{lazar_flexible_2009,balau_one_2011,hermans_horizon-1_2013} or  added as constraints~\cite{wu_control_2019,zeng_enhancing_2021,ma_feasibility_2021}. However,  {those} approaches result again in nonlinear  {and} nonconvex optimization problems. }

 {In summary, standard MPC provides stability and safety guarantees but its computational burden is prohibitive for applications with severely constrained computational resource. This holds despite real-time schemes  {such as} RTI which reduce the computational effort but generally lack runtime guarantees  {and often} sacrifice optimality  {or} feasibility. Conversely, CLF/CBF methods offer  {QP-based feedback laws with small computational footprints} and potential runtime guarantees, yet lack comparable theoretical guarantees for stability and safety under nonlinear dynamics.
This paper closes this gap by fusing the advantages of the convex QP formulation  {of CLF/CBF methods} with the stability and safety guarantees of quasi-infinite horizon MPC, resulting in a computationally lightweight scheme, which we refer to as {\em infinitesimal-horizon} MPC ($\partial$MPC).} The main contributions of this paper are:
\paragraph{Real-time feasibility}  An efficient QP formulation for the infinitesimal-horizon MPC implies safety and {asymptotic stability} guarantees in real time.
\paragraph{Theoretical guarantees}  {The  {infinitesimal-horizon} MPC  {scheme}  {is proven to be} recursive {ly feasible} and asymptotic {ally stable} under mild assumptions.}
\paragraph{Practical synthesis} A nonlinear sum-of-squares optimization method synthesizes approximations that satisfy the sufficient conditions for safety and asymptotic stability.

The proposed approach ensures safety of a real-time feasible approximate MPC {\em by design}.  {Whereas our feedback law is based on an instantaneous optimization similar to CLF/CBF approaches \cite{ames_control_2017, schneeberger_advanced_2026}, we neither require nor assume continuous feedback, allowing for MPC-inspired cost functions and constraints.} 
 {We demonstrate that}, in comparison to state-of-the-art methods for a nonlinear constrained satellite attitude control problem,  {our approach}  {significantly reduces}  {the} computation time while being competitive in performance.
Specifically, we demonstrate the effectiveness of our approach by enabling a spacecraft to perform challenging large-angle, three-axis maneuvers  {under nonconvex constraints that protect} sensitive instruments like telescopes. 

The remainder of this paper is as follows: Section~\ref{sec: ProblemStatement} provides the problem statement and mathematical background, Section~\ref{sec: Methodology} provides the theoretical results of this paper and the QP formulation. Section~\ref{sec: SOSApproxAndQP} provides a sum-of-squares (SOS) program to synthesize the  {needed ingredients}.  {Our} approach is demonstrated and compared to other methods in Section~\ref{sec: NumericalResults}.

\section{Problem Statement}
\label{sec: ProblemStatement}
We consider the continuous-time nonlinear dynamic system
\begin{align}
    \label{eq:system}
    \dot x(t) = f(x(t), u(t))
\end{align}
for (almost) all $t \geq 0$ with states $x(t) \in \mathbb R^n$ and inputs $u(t) \in \mathcal U$, where $\mathcal U \subset \mathbb R^m$ denotes the viable inputs. 
We say that \eqref{eq:system} is {\em control-affine} if $f(x,u)$ is affine in $u$ for all $x \in \mathbb R^n$.

Throughout the paper we assume that $f: \mathbb R^n \times \mathcal U \to \mathbb R^n$ is Lipschitz continuous, satisfies $f(0,0) = 0$, and $\mathcal U$ is a compact set with $0 \in \interior \mathcal U$. For any $t_1 \geq t_0 \geq 0$, we denote by $\mathscr U_{[t_0, t_1]}$ the set of Lebesgue-measurable functions $\uu: [t_0, t_1] \to \mathcal U$.  { {For} any $\uu\in\Ufcnset{t_0,t_1}$ let
$\xi(t,\uu,t_0,x_0)$ denote the solution of the initial value problem of
\eqref{eq:system} on $[t_0,t_1]$ with $x(t_0)=x_0$ under the input
$u(t)=\uu(t)$.}

 { {The} feedback laws considered in this paper
are obtained pointwise as minimizers of an optimization problem and hence, the closed-loop initial value problem may
admit no solution in the classical or Carath\'eodory sense. We therefore adopt
the notion of Krasovskii {\cite{ceragioli1999discontinuous,cortes2008discontinuous}}:  {Take}
$\mathcal D\subseteq\mathbb R^n$ and let $F:\mathcal D\to\mathbb R^n$ be
locally bounded. For  {any} $x\in\mathcal D$, let
\begin{align}
    \label{eq:krasovskii}
    \mathscr K[F](x) := \bigcap_{\varepsilon>0}
        \overline{\operatorname{co}}\, F\bigl(B_\varepsilon(x)\cap\mathcal D\bigr),
\end{align}
where $B_\varepsilon(x)$ is the open ball of radius $\varepsilon$  {centered at} $x$;
the set $\mathscr K[F](x)$ is then nonempty, compact, and convex, and
$\mathscr K[F]$ is upper semicontinuous on $\mathcal D$. Given
$\kappa:\mathcal D\to\Uset$,  {locally bounded by compactness of $\Uset$}, a \emph{closed-loop solution} of \eqref{eq:system} under $\kappa$ is
an absolutely continuous function $\mathbf x$ with values in $\mathcal D$ such
that
\begin{align}
    \label{eq:closed-loop-inclusion}
    \dot{\mathbf x}(t)\in\mathscr K[f(\cdot,\kappa(\cdot))](\mathbf x(t))
    \quad\text{for almost all } t .
\end{align}
We write $\phi(t,\kappa,t_0,x_0)$ for the value at time $t$ of a {maximal}
closed-loop solution with $\mathbf x(t_0)=x_0$; every statement involving
$\phi$ is asserted for \emph{every} such solution. Note that no measurability of $\kappa$ is required, since \eqref{eq:krasovskii}
discards no sets of measure zero; in particular $F(x)\in\mathscr K[F](x)$ for
every $x\in\mathcal D$. If $\kappa$ is continuous at $x$, then
$\mathscr K[f(\cdot,\kappa(\cdot))](x)=\{f(x,\kappa(x))\}$ and
\eqref{eq:closed-loop-inclusion} reduces to the usual closed-loop dynamics.}

\subsection{Model Predictive Control}
Let $\CS, \TS \subset \mathbb R^n$ be the closed state and terminal constraint sets, respectively (with $0\in\textrm{int}~\TS$), and $T > 0$ be the prediction horizon. The  MPC optimal value function for the state $x_k\in\CS$ is given by
\begin{subequations}
    \label{eq:mpc-problem}
\begin{align}
    V(x_k) = {} &\min_{\mathbf u \in \mathscr U_{[0, T]}} \mathscr J_{[0, T]}(\xi(\cdot, \mathbf u, 0, x_k), \mathbf u) \\
    &\text{subject to $\xi(T, \mathbf u, 0, x_k) \in \TS$} \\
    &\text{and $\forall t \in [0, T], \, \xi(t, \mathbf u, 0, x_k) \in \CS$}
\end{align}
with cost functional
\begin{align}
    \label{eq:mpc-cost}
    \mathscr J_{[0, T]}: (\mathbf x, \mathbf u) \mapsto F(\mathbf x(T)) + \int_0^T L(\mathbf x(t), \mathbf u(t)) \mathrm d t,
\end{align}
\end{subequations}
{continuously differentiable}  {\em terminal cost} $F:  {\TS} \to \mathbb R_{\geq 0}$, and  {continuous  {\em stage cost}} $L: {\CS} \times \mathcal U \to \mathbb R_{\geq 0}$ (both positive definite w.r.t. their respective origins).\footnote{A common choice is $L: (x,u) \mapsto x^\top Q x + u^\top R u$ and $F: x \mapsto x^\top S x$ with positive definite matrices $Q, S \in \mathbb R^{n \times n}$ and $R \in \mathbb R^{m \times m}$.}

 {The analysis of MPC schemes is contingent on the feasibility of the optimization problem in \eqref{eq:mpc-problem}. For the prediction horizon $T$ and an initial condition $x \in \CS$, the set of {\em feasible solutions} to \eqref{eq:mpc-problem} is
\begin{multline}
    \label{eq:mpc-feasible}
    \mathscr F(x) = \{ \uu\in\Ufcnset{0,T} 
    ~|~ \xi(T,\uu,0,x)\in\TS \\ \text{and\;}\forall t\in[0,T],\,\xi(t,\uu,0,x)\in\CS\}
\end{multline}
and $\Feas = \{ x \in \CS ~|~ \mathscr F(x) \neq \varnothing \}$ is called the {{\em set of feasible initial conditions}}.
The following result associates feasibility of \eqref{eq:mpc-problem} to the properties of the nonlinear system \eqref{eq:system} relative to the state and terminal constraints.}

\begin{lemma}[\cite{Cunis2021aut}]
     {The set of feasible initial conditions $\Feas$ is equal to the set
    \begin{multline}
        \label{eq:sys-reachavoid}
        \RA{t_0,t_1}{\CS,\TS} = \{x_0\in\CS ~|~ \exists \uu\in\Ufcnset{t_0,t_1}, \forall t\in[t_0,t_1], \\ \text{$\xi(t,\uu,t_0,x_0)\in\CS$ and $\xi(t_1,\uu,t_0,x_0)\in\TS$} \}.
    \end{multline}
    for $t_0=0$ and $t_1=T$.}
\end{lemma}

 {The set $\RA{t_0,t_1}{\CS,\TS}$ is referred to as reach-avoid set or, in the context of MPC, as controllability set of the system \eqref{eq:system}.}
 { {Being} feasible for any given $x_k \in \Feas$,  {we assume that \eqref{eq:mpc-problem}} admits an optimal solution $\mathbf u_{x_k}^\star \in \mathscr U_{[0, T]}$.}
Given a sequence $\{ t_k \}_{k \in \mathbb N_0} \subset \mathbb R_{\geq 0}$ of time instances with $t_k < t_{k+1} \leq t_k + T$ for all $k \in \mathbb N_0$ and an initial condition $x_0 \in \CS$,
we define%
\begin{align}
    \label{eq:mpc-feedback}
    u_\text{MPC}: t\mapsto \mathbf u_{x_k}^\star(t-t_k), \quad t \in [t_k, t_{k+1}),
\end{align}
where $\mathbf u_{x_k}^\star$ is the optimal solution for $x_k$ (feasibility provided) and $x_{k+1} = \xi(t_{k+1}, u_\text{MPC}(\cdot), t_k, x_k)$ is the next sample, for all $k \in \mathbb N_0$. This results in a {\em receding-horizon} MPC scheme.  {When considering closed-loop operation for all $t\ge 0$, we additionally assume $t_k\to\infty$.}

\subsection{Control Barrier and Lyapunov Functions}


    {A continuously differentiable function} $\clf :\CS\rightarrow \mathbb{R}_{\geq 0}$ {is a control Lyapunov function \cite{Sontag1999} of \eqref{eq:system} } if  $\clf$ is \emph{proper}, \ie, all of its sublevel sets are compact;
      $\clf$ is positive on $\CS \setminus \{0\}$ with $\clf[0] = 0$; and
       for all $x\in\CS\setminus\{0\}$, there exists some $u\in\Uset$ such that 
            \begin{align}
                \label{eq:CLFcondition}
                \dot V(x,u) = \ip{\nabla V(x),f(x,u)}<0 {.}
            \end{align}
%
%
%
%
    {Similarly, a continuously differentiable function} $\cbf :\CS\rightarrow\mathbb{R}$ {is a control barrier function \cite{ames_control_2019} of \eqref{eq:system}} if there exists an extended class $\mathcal{K}$\footnote{ {An \emph{extended class} $\mathcal{K}$ function is any continuous strictly increasing function $\alpha:\mathbb{R}\to \mathbb{R}$ that satisfies $\alpha(0)=0$.}} function $\alpha$ such that for every $x\in\CS$, there exists a $u\in\Uset$ such that 
    \begin{align}
        \label{eq:CBFcondtion}
        \dot \cbf(x,u) = \ip{\nabla \cbf[x],f(x,u)}\leq \alpha(-\cbf[x]) {.}
    \end{align}
%
The existence of a CBF guarantees that  {for any $x_0 \in \mathbb R^n$ with $h(x_0) \leq 0$ and $T > 0$ there exists $\mathbf u \in \mathscr U_{[0, T]}$ such that $\xi(\cdot, \mathbf u, 0, x_0)$ remains in} the sublevel set $\{x\in\CS~|~h(x)  {\leq} 0\}$. 

Oftentimes,  {CBF and CLF are used simultaneously} to guarantee safety and  {stabilizability}  {(cf.~\cite{ames_control_2019})}. In those cases, it is not enough to have {a pair of distinct CLF and CBF; rather, they must be \emph{compatible}, \ie, for every $x\in\CS$ there must exist a $u\in\Uset$ such that \eqref{eq:CLFcondition} and \eqref{eq:CBFcondtion} hold simultaneously.  {Note that an  {relaxation} of this problem is to employ a slack variable that  {balances} stabilizability and safety \cite{ames_control_2019}, but this does not resolve the issue of potential incompatibility.}}
 {In the next section we will show that the MPC \eqref{eq:mpc-problem} can be efficiently solved by  compatible CLF/CBF pairs for \eqref{eq:system}.} 

\subsection{ {Relaxed} Dynamic Programming Principle}
\label{subsec: relaxedDP}
In MPC, the objective of minimizing the integral of a cost over a finite horizon $T$  {must be balanced against stabilization of the closed-loop dynamics and satisfaction of} state constraints.  {Otherwise}, the {finite-horizon} cost of operation is not a good measure for an asymptotically well-behaved dynamic  {(see, e.g., \cite[Example~7.1]{gruenePannek2017})}.  {In this work, we} ensure asymptotic stability and constraint satisfaction  {by subjecting} the terminal cost $F$ and the terminal set $\TS$ in \eqref{eq:mpc-problem} to  {a} stabilizing terminal condition  {following} \cite[Assumption~5.9]{gruenePannek2017}. 
 {As we are going to refer to this condition for different pairs of costs $P$ and sets $\Omega$, we introduce the following definition.}
 {Recall that the {tangent cone} of a set $\mathcal A \subset \mathbb R^n$ at $\bar x \in \mathcal A$ is the set $\mathcal T_{\bar x} \mathcal A$ of all vectors $z \in \mathbb R^n$ for which there exists sequences $x_k \in \mathcal A$ and  {$h_k > 0$} with $x_k \to \bar x$,  {$h_k \to 0$}, and $(x_k - \bar x)/h_k \to z$.}

\begin{definition}
    \label{def:terminal}
    {A pair $(P, \Omega)$ with $0 \in \interior \Omega$, $\Omega\subseteq\CS$, and  {differentiable} $P: \Omega \to \mathbb R$ satisfies the {\em stabilizing terminal condition} for the system \eqref{eq:system} if} 
    there exists a  {continuous} state feedback $\kappa: \Omega \to \mathcal U$ such that
    \begin{subequations} 
        \label{eq:terminal}
    \begin{align}
        \label{eq:terminal-cost}
        \langle\nabla P(x) ,f(x, \kappa(x))\rangle &\leq -L(x, \kappa(x)) \\
        \label{eq:terminal-set}
        f(x,\kappa(x)) &\in\Tbdl{x}\Omega
    \end{align}
    \end{subequations}
    for all $x \in \Omega$, where $\Tbdl{x}\Omega$ is the tangent cone of $\Omega$ at $x\in\Omega$.
\end{definition}


 {A well-known  {result} is that, if $(F, \TS)$ satisfy the  {stabilizing terminal condition} in Definition~\ref{def:terminal}, then \eqref{eq:mpc-problem} remains feasible for all $x_k$, $k \in \mathbb N_0$, with $x_0 \in \mathcal R$ and the closed-loop of \eqref{eq:system} and \eqref{eq:mpc-feedback} is asymptotically stable on $\mathcal R$ \cite{magni_stabilizing_2004, Mayne2000}.}
 {The mechanism behind stabilizing MPC via a terminal condition is formalized in the following lemma.}

\begin{lemma}
    \label{lem:infinite-ddp}
    Let $(F, \TS)$ satisfy the stabilizing terminal condition of Definition~\ref{def:terminal}; then
    \begin{subequations}
        \label{eq:infinite-ddp}
    \begin{align}
        V(x_0) \geq {} &\inf_{\mathbf u \in \mathscr U_{[t_0, t_1]}} \mathscr V_{[t_0, t_1]}(\xi(\cdot, \mathbf u, t_0, x_0), \mathbf u) \\ 
        &\text{subject to $\forall t \in [t_0, t_1], \, \xi(t, \mathbf u, t_0, x_0) \in \Feas$}
        \label{eq:infinite-ddp-reach}
    \end{align}
    \end{subequations}
    with
    \begin{align*}
        \mathscr V_{[t_0, t_1]}: (\mathbf x, \mathbf u) \mapsto V(\mathbf x(t_1)) + \int_{t_0}^{t_1} L(\mathbf x(t), \mathbf u(t)) \mathrm d t
    \end{align*}
    for all $t_1 > t_0 \geq 0$ and $x_0 \in \mathcal R$.
\end{lemma}

 {The proof is standard (cf.~\cite{chen1998}). Since \eqref{eq:infinite-ddp} holds for arbitrarily short $[t_0,t_1]$, and if $V$ and $\Feas$ are a priori known, solving \eqref{eq:infinite-ddp} instead of \eqref{eq:mpc-problem} could reduce computational effort. However, computing these exactly is difficult or impossible, and $V$ is generally nonsmooth. To address this, we next consider sufficiently regular approximations of $V$ and $\Feas$, show how a control law can be efficiently computed for them, and that this still guarantees good properties for \eqref{eq:mpc-problem}.}

\begin{remark}
    {If $(F, \TS)$ satisfies the stabilizing terminal condition,}
    then  {Lemma~\ref{lem:infinite-ddp}} implies that the terminal penalty $F$ overapproximates the  {optimal running} cost  {$\int L(\mathbf x(t), \mathbf u(t)) \mathrm d t$} on the terminal set $\TS$,  {where the local state feedback $\kappa$ in Definition~\ref{def:terminal} provides an suboptimal solution to \eqref{eq:infinite-ddp} that keeps the system \eqref{eq:system} within $\Omega$ (a fortiori $\mathcal X$).
    Consequently, $V$ and $\mathcal R$ become approximations of the {\em infinite-horizon} optimal cost and controllable set, respectively.}
\end{remark}



\section{Methodology}
\label{sec: Methodology}

 {In this section, we will show that instead of solving the MPC problem with the exact value function $V$ and feasible set $\Feas$, one can instead solve a more tractable optimization problem using  {a} smooth positive-definite $\Vap$ and $\Rap=\{x\in\mathbb{R}^n~|~\hap(x)\leq0\}$ for some continuously differentiable 
{proper} $\hap:\mathbb{R}^n\to\mathbb{R}$ such that $\Rap\subseteq\CS$ as long as the pair $(\Vap,\Rap)$ satisfies the stabilizing terminal condition. More precisely, given $\alpha$ of  {extended} class $\mathcal{K}$, we consider the infinitesimal MPC ($\partial$MPC) problem defined as 
\begin{subequations}
    \label{eq:mpc-infinitesimal}
    \begin{align}
        &\underset{u\in\Uset}{\text{minimize}}\quad L(x,u)+\ip{\nabla \Vap(x),f(x,u)}\\
        &\text{subject to}~~\ip{\nabla \hat h(x), f(x,u)} \leq \alpha( -\hat h{(x)} ),\label{eq:mpc-infinitesimal-constraint}
    \end{align}
\end{subequations}
whose optimal solution set is denoted $ {\mathcal S_\partial}(x)$.  {Our main result is as follows.}}

\begin{theorem}
    \label{thm:main}
     {Let $\Rap:=\{x\in\mathbb R^n\mid\hap(x)\le0\}$  {be compact} for
    some continuously differentiable  $\hap$ that satisfies
    $\nabla\hap(x)\neq0$ whenever $\hap(x)=0$, and let $\Vap$ be continuously
    differentiable and positive definite on a neighbourhood of $\Rap$. Assume
    the pair $(\Vap,\Rap)$ satisfies the stabilizing terminal condition of
    Definition~\ref{def:terminal}. Then there exists an extended class
    $\mathcal K$ function $\alpha$ such that the
    infinitesimal MPC problem \eqref{eq:mpc-infinitesimal}
     {satisfies}:}
    \begin{enumerate}
        \item  {for every $x\in\Rap$, \eqref{eq:mpc-infinitesimal}
            is feasible and 
            $ {\mathcal S_\partial}(x)\neq\varnothing$.}
    \end{enumerate}
     {Let now $\uuinf:\Rap\to\Uset$ be any selection of the optimal
    solution map, \ie, $\uuinf(x)\in {\mathcal S_\partial}(x)$ for all $x\in\Rap$;
    then moreover:}
    \begin{enumerate}
        \setcounter{enumi}{1}
        \item  {for every $x_0\in\Rap$  {and $t_0 \geq 0$}, at least one closed-loop
            solution under $\uuinf$  {in the sense of \eqref{eq:closed-loop-inclusion} } exists,
            and every maximal  {solution} satisfies
            $\phi(t,\uuinf,t_0,x_0)\in\Rap$ for all $t\ge t_0$;}
        \item  {the origin is Lyapunov stable  {for
            the closed loop} and every maximal
            closed-loop solution  {satisfies $\phi(t, \uuinf,t_0,x_0) \to 0$ as $t \to \infty$ if $x_0 \in \Rap$.}}
    \end{enumerate}
\end{theorem}

 {The assumptions warrant some remarks. $\Rap$ is compact if $\CS$ is compact, since $\hap$ is continuous and $\Rap \subseteq \CS$ by the stabilizing terminal condition. The assumption of nonvanishing gradients $\nabla \hap$ on its zero levelset implies that the latter is exactly the boundary of $\Rap$. Our assumption here is similar to those of \cite{ames_control_2017, schneeberger_advanced_2026} and is usually satisfied by polynomial base functions. In the subsequent section, we will propose a nonconvex SOS optimization framework to establish polynomial approximations $\Vap$ and $\hap$ subject to the stabilizing terminal condition.}

 {Under the assertions of Theorem~\ref{thm:main},  {which we will prove through multiple lemmas,} the feedback defined by \eqref{eq:mpc-infinitesimal} is recursively feasible on $\Rap \subseteq \CS$ and (strongly) asymptotically stabilized \eqref{eq:system} with region of attraction $\Rap$.}
 {
 {We start by looking}
at how the dynamic programming principle for the approximation pair $(\Rap,\Vap)$ can be used to motivate the $\partial$MPC problem, and how it relates to notions of control barrier functions.}


\subsection{ {Infinitesimal-horizon MPC}}
 {Lemma~\ref{lem:infinite-ddp} holds for every horizon $t_1 > t_0$, and
letting that horizon shrink is what motivates the form of
\eqref{eq:mpc-infinitesimal}.  {In particular, it is easy to see that \eqref{eq:infinite-ddp} is satisfied with $(V, \mathcal R)$ replaced by $(\Vap, \Rap)$ provided that the latter satisfy the stabilizing terminal condition.} Write  {now} $t_1 = t_0 + h$  {with $h > 0$} and fix $x_0\in\Rap$;
since $\Vap$ is continuously differentiable and $L$ is continuous, any input $\uu$  {and $\xx = \xi(\cdot,\uu,t_0,x_0)$}
satisfy
\begin{multline*}
    \Vap(\xx(t_0+h)) + \int_{t_0}^{t_0+h} L(\xx(t),\uu(t))\dd t - \Vap(x_0)
        \\= \ip{\nabla\Vap(x_0), f(x_0,\uu(t_0))} h + L(x_0,\uu(t_0)) h + o(h).
\end{multline*}
 {Moreover,} the state constraint $\xi(t,\uu,t_0,x_0)\in\Rap$  becomes, after subtracting $x_0$ and
dividing by $h$, a condition on the velocity $f(x_0,\uu(t_0))$ relative to
$\Tbdl{x_0}\Rap$. Substituting  {for $V$ in \eqref{eq:infinite-ddp}},
canceling $\Vap(x_0)$, and dividing by $h$, the relaxed dynamic programming
principle collapses as  {$h \to 0$} to the pointwise statement
\begin{align*}
    \min_u\; L(x_0,u) + \ip{\nabla\Vap(x_0), f(x_0,u)}
    \;\le\; 0,
\end{align*}
the minimum being taken over those $u\in\Uset$  {with $f(x_0, u) \in \Tbdl{x_0}\Rap$.} 
The $\partial$MPC problem \eqref{eq:mpc-infinitesimal} is exactly
this limit, with the tangent-cone cone traded for the barrier.}


\begin{lemma}
    \label{lem:cone-to-cbf}
     {Let $\Rap=\{x\in\mathbb{R}^n\mid\hap(x)\le0\}$ be compact with
    $\hap$ continuously differentiable, and let $\kappa:\Rap\to\Uset$ be
    continuous with $f(x,\kappa(x))\in\Tbdl{x}\Rap$ for all $x\in\Rap$. Then
    there exists an extended class $\mathcal K$ function $\alpha$ with}
    \begin{align}
        \label{eq:terminal-classK}
        \ip{\nabla\hap(x), f(x,\kappa(x))} \le \alpha(-\hap(x))
    \end{align}
     {for all $x \in \Rap$.}
\end{lemma}
\begin{proof}
    {Let $\psi:=\ip{\nabla\hap(\cdot),f(\cdot,\kappa(\cdot))}$,
    continuous on the compact set $\Rap$. We will divide the proof into three steps.
    \emph{(i)} $\psi(x)\le0$ for $x\in Z:=\{x\in\Rap\mid\hap(x)=0\}$. Indeed, if
    $x\in\partial\Rap$ then $\Tbdl{x}\Rap\subseteq\{v\mid\ip{\nabla\hap(x),v}\le0\}$,
    since any $v\in\Tbdl{x}\Rap$ is a limit $(x_k-x)/h_k$ with $x_k\in\Rap$ and
    $\hap(x_k)\le0=\hap(x)$, so that
    $\ip{\nabla\hap(x),x_k-x}\le o(\|x_k-x\|)$; if instead $x\in\interior\Rap$
    then $\hap$ attains an interior maximum at $x$, so $\nabla\hap(x)=0$.
    \emph{(ii)} Let $s_{\max}:=\max_{\Rap}(-\hap)$ and, for
    $s\in(0,s_{\max}]$, let $\Psi(s):=\max\{\psi(x)\mid x\in\Rap,\ -\hap(x)=s\}$
    (and $\Psi(s):=0$ if the level set is empty), which is attained and bounded
    above by $\max_{\Rap}\psi<\infty$. Then $\limsup_{s\to0}\Psi(s)\le0$:
    for $s_k\to0$ and $x_k$ attaining $\Psi(s_k)$, compactness gives a
    subsequence $x_k\to\bar x$ with $\hap(\bar x)=0$, whence
    $\Psi(s_k)=\psi(x_k)\to\psi(\bar x)\le0$ by (i).
    \emph{(iii)} Set $\beta(0):=0$ and
    $\beta(s):=\max\{0,\sup_{0<\sigma\le s}\Psi(\sigma)\}$, which is
    non-decreasing and bounded, and such that there exists $\delta>0$ for which $\beta$ is continuous in $[0,\delta)$ by \emph{(ii)}. Define
    $\alpha(s):=\int_1^2\beta(s\tau)\,\dd\tau+\varepsilon s$ for $s\ge0$ and
    $\alpha(s):=-\alpha(-s)$ for $s<0$, with  {some} $\varepsilon>0$. Since $\beta$
    is non-decreasing, $s\mapsto\beta(s\tau)$ is non-decreasing for each
    $\tau\in[1,2]$, so $\alpha$ is non-decreasing, and the term $\varepsilon s$
    makes it strictly increasing; $\alpha$ is continuous with $\alpha(0)=0$
    because $\beta(0^+)=0$; and $\alpha(s)\ge\int_1^2\beta(s)\dd\tau=\beta(s)$.
    Therefore, if $\hap(x)=0$ then $\psi(x)\le0=\alpha(0)=\alpha(-\hap(x))$ by
    \emph{(i)}, and if $\hap(x)<0$ then $s:=-\hap(x)\in(0,s_{\max}]$ and
    $\psi(x)\le\Psi(s)\le\beta(s)\le\alpha(s)$.}
\end{proof}

 {We now prove the first claim of Theorem~\ref{thm:main}, the feasibility of \eqref{eq:mpc-infinitesimal} on $\Rap$.}

\begin{lemma}
    \label{lem: feasInfMPC}
     {Let $\Rap=\{x\in\mathbb{R}^n\mid\hap(x)\le0\}$  {be compact} for some continuously
    differentiable  $\hap$, let $\Vap:\Rap\to\mathbb R$ be
    differentiable, and assume the pair $(\Vap,\Rap)$ satisfies the stabilizing
    terminal condition of Definition~\ref{def:terminal}. Then there exists an
    extended class $\mathcal K$ function $\alpha$ such that, for every
    $x\in\Rap$, the $\partial$MPC problem \eqref{eq:mpc-infinitesimal} is
    feasible and attains its minimum; that is,
    $ {\mathcal S_\partial}(x)\neq\varnothing$.}
\end{lemma}
\begin{proof}
     {
    Let $\kappa:\Rap\to\Uset$ be a continuous state feedback satisfying
    \eqref{eq:terminal} for all $x\in\Rap$, which exists by
    Definition~\ref{def:terminal}, and let $\alpha$ be the extended class
    $\mathcal K$ function that Lemma~\ref{lem:cone-to-cbf} associates with
    $\kappa$. 
    Fix
    $x\in\Rap$. Lemma~\ref{lem:cone-to-cbf} gives
    $\ip{\nabla\hap(x),f(x,\kappa(x))}\le\alpha(-\hap(x))$, so $u=\kappa(x)$
    satisfies \eqref{eq:mpc-infinitesimal-constraint} and the feasible set}
    \begin{align*}
        { \Gamma(x) := \bigl\{u\in\Uset \bigm|
            \ip{\nabla\hap(x),f(x,u)}\le\alpha(-\hap(x))\bigr\}}
    \end{align*}
     {is nonempty. It is closed, being the intersection of $\Uset$
    with the preimage of a closed half-line under a map continuous in $u$, and
    hence compact because $\Uset$ is compact. Since the objective of
    \eqref{eq:mpc-infinitesimal} is continuous in $u$, it attains its minimum
    on $\Gamma(x)$, so $ {\mathcal S_\partial}(x)\neq\varnothing$.}
\end{proof}

 {Henceforth, let $\alpha$ be the extended class $\mathcal K$ function associated with the feedback $\kappa$ from Definition~\ref{def:terminal} in the sense of Lemmas~\ref{lem:cone-to-cbf} and \ref{lem: feasInfMPC}.}

\subsection{ {Recursive Feasibility}}
 {In this subsection we show  {that} the closed-loop under the solution of \eqref{eq:mpc-infinitesimal} is  {forward complete and} recursively feasible.}

\begin{lemma}
\label{lem: recFeasInfMPC}
 {Let $\Rap=\{x\in\mathbb{R}^n\mid\hap(x)\le0\}$  {be compact} for some
continuously differentiable  $\hap$ satisfying $\nabla\hap(x)\neq0$
whenever $\hap(x)=0$, let $\Vap:\Rap\to\mathbb R$ be differentiable, and assume
the pair $(\Vap,\Rap)$ satisfies the stabilizing terminal condition of
Definition~\ref{def:terminal}. Let $\uuinf:\Rap\to\Uset$  $\uuinf(x)\in {\mathcal S_\partial}(x)$ for all
$x\in\Rap$. Then, for every
$t_0\ge0$ and every $x_0\in\Rap$,}
\begin{enumerate}
    \item  {at least one closed-loop solution
        $\phi(\cdot,\uuinf,t_0,x_0)$ exists}; 
    \item  {every maximal closed-loop solution is defined on
        $[t_0,\infty)$ and satisfies $\phi(t,\uuinf,t_0,x_0)\in\Rap$ for all
        $t\ge t_0$.}
\end{enumerate}
 {In particular, \eqref{eq:mpc-infinitesimal} is feasible at
$\phi(t,\uuinf,t_0,x_0)$ for every $t\ge t_0$, so the $\partial$MPC scheme is
recursively feasible.}
\end{lemma}

\begin{proof}
 {Let $F:=f(\cdot,\uuinf(\cdot)):\Rap\to\mathbb R^n$. Since $\Rap$ is compact, $M:=\max_{\Rap\times\Uset}\|f\|$ is finite and
$\|F\|\le M$ on $\Rap$; in particular, $F$ is bounded and $\mathscr K[F]$ has
nonempty, compact, and convex values and is upper semicontinuous on $\Rap$.
We first record an estimate that will be used twice. Let $c:\Rap\to\mathbb R^n$
and $\rho:\Rap\to\mathbb R$ be continuous and suppose
$\ip{c(y),F(y)}\le\rho(y)$ for all $y\in\Rap$. Fix $x\in\Rap$, let
$\varepsilon>0$, and let $y\in B_\varepsilon(x)\cap\Rap$; then
\begin{align*}
    \ip{c(x),F(y)}
    &= \ip{c(y),F(y)} + \ip{c(x)-c(y),F(y)} \\
    &\le \rho(y) + \|c(x)-c(y)\|\,M \;\le\; \rho(x)+\omega(\varepsilon),
\end{align*}
where $\omega(\varepsilon)\to0$ as $\varepsilon\to0$ because $c$ and
$\rho$ are uniformly continuous on the compact set $\Rap$. As
$v\mapsto\ip{c(x),v}$ is linear and continuous, this bound extends to
$\overline{\operatorname{co}}\,F(B_\varepsilon(x)\cap\Rap)$, and letting
$\varepsilon\to0$ yields
\begin{align}
    \label{eq:inheritance}
    \ip{c(x),v}\le\rho(x)\qquad\text{for all } v\in\mathscr K[F](x).
\end{align}
Since $\uuinf(y)\in {\mathcal S_\partial}(y)$ is feasible for
\eqref{eq:mpc-infinitesimal} at $y$, applying \eqref{eq:inheritance} with
$c=\nabla\hap$ and $\rho=\alpha\circ(-\hap)$ gives
\begin{align}
    \label{eq:cbf-inherited}
    \ip{\nabla\hap(x),v}\le\alpha(-\hap(x))
    \qquad\text{for all } v\in\mathscr K[F](x).
\end{align}
Now let $x\in\Rap$ with $\hap(x)=0$. Then $\alpha(-\hap(x))=\alpha(0)=0$, so
\eqref{eq:cbf-inherited} gives $\ip{\nabla\hap(x),v}\le0$ for every
$v\in\mathscr K[F](x)$; since $\nabla\hap(x)\neq0$ by hypothesis,
$\Tbdl{x}\Rap=\{v\mid\ip{\nabla\hap(x),v}\le0\}$ and therefore
$\mathscr K[F](x)\subseteq\Tbdl{x}\Rap$. At any $x$ with $\hap(x)<0$ we have
$\Tbdl{x}\Rap=\mathbb R^n$. In either case
$\mathscr K[F](x)\cap\Tbdl{x}\Rap\neq\varnothing$, so $\Rap$ is viable for the
differential inclusion \eqref{eq:closed-loop-inclusion} and, by the viability
theorem for upper semicontinuous inclusions with nonempty compact convex values
on a closed set, for every $x_0\in\Rap$ there is a
closed-loop solution with $\mathbf x(t_0)=x_0$, which proves the first claim.
Finally, let $\mathbf x$ be any maximal closed-loop solution, defined on some
$[t_0,t_{\max})$ and taking values in $\Rap$, and suppose $t_{\max}<\infty$.
Then $\|\dot{\mathbf x}\|\le M$ almost everywhere, so $\mathbf x$ is
$M$-Lipschitz and hence uniformly continuous, its limit
$\bar x:=\lim_{t\to t_{\max}}\mathbf x(t)$ exists and lies in the closed
set $\Rap$, and setting $\mathbf x(t_{\max}):=\bar x$ leaves $\mathbf x$ a
solution. Applying the viability theorem at $\bar x$ and concatenating extends
$\mathbf x$ beyond $t_{\max}$, contradicting maximality. Hence
$t_{\max}=\infty$, and $\mathbf x(t)\in\Rap$ for all $t\ge t_0$.}
\end{proof}


 {Since $\uuinf$ is defined on $\Rap$ only, closed-loop solutions in
the sense of \eqref{eq:closed-loop-inclusion} take their values in $\Rap$ by
construction. The content of Lemma~\ref{lem: recFeasInfMPC} is therefore that
the closed loop never  {leaves the domain of $\uuinf$ nor} 
escapes in finite
time. Together with $\Rap\subseteq\CS$ this yields constraint satisfaction for
all $t\ge t_0$.}

\subsection{ {Asymptotic Stability}}
 {Finally, we tackle the last claim of Theorem \ref{thm:main}, showing that the closed loop is not only stable, but  {the origin is (strongly) asymptotically stable with region of attraction $\Rap$. This is because the origin} is the  {(unique)} minimizer of $L$  {which} is assumed positive definite.}
\begin{lemma}
\label{lem: optValDecInfMPC}
 {Let the hypotheses of Lemma~\ref{lem: recFeasInfMPC} hold, with
$\Vap$ in addition continuously differentiable and positive definite. Then the origin is Lyapunov stable  {for the closed loop \eqref{eq:closed-loop-inclusion}} and every
maximal closed-loop solution 
with $x_0\in\Rap$
satisfies $\phi(t,\uuinf,t_0,x_0)\to 0$ as $t\to\infty$.}
\end{lemma}

\begin{proof}
 {Write $F(x):=f(x,\uuinf(x))$ and
$W(x):=\min_{u\in\Uset}L(x,u)$. As $L$ is continuous and $\Uset$ is
compact, $W$ is continuous, and it
is positive definite because $L$ is. Fix $x\in\Rap$. By
Definition~\ref{def:terminal} there is some $u\in\Uset$ feasible for
\eqref{eq:mpc-infinitesimal} at $x$ with
$\ip{\nabla\Vap(x),f(x,u)}+L(x,u)\le 0$, so the optimal value of
\eqref{eq:mpc-infinitesimal} is nonpositive; since $\uuinf(x)$ attains it,
\begin{align}
    \label{eq:V-decrease-pointwise}
    \ip{\nabla\Vap(x),F(x)}\le-L(x,\uuinf(x))\le-W(x)
\end{align}
 {for all $x\in\Rap$.}
Both $\nabla\Vap$ and $W$ are continuous, so the estimate
\eqref{eq:inheritance} established in the proof of
Lemma~\ref{lem: recFeasInfMPC}, applied with $c=\nabla\Vap$ and $\rho=-W$,
transfers \eqref{eq:V-decrease-pointwise} to the regularisation,
\begin{align}
    \label{eq:V-decrease-inherited}
    \ip{\nabla\Vap(x),v}\le-W(x)
    \qquad\text{for all } v\in\mathscr K[F](x)
\end{align}
 {for all $x\in\Rap$.}}
 {Then, to verify that the origin is an equilibrium of the closed loop, notice that
$0\in\operatorname{int}\Rap$ and $\Vap$ attains its minimum over $\Rap$ there,
so $\nabla\Vap(0)=0$. Since $\Rap$ and $\Uset$ are compact,
$M:=\max_{x\in\Rap,u\in\Uset}\|f(x,u)\|<\infty$, and
\eqref{eq:V-decrease-pointwise} gives
$L(y,\uuinf(y))\le\|\nabla\Vap(y)\|\,M\to 0$ as $y\to0$ in $\Rap$. Positive
definiteness of $L$ then forces $\uuinf(y)\to0$, whence $F(y)\to f(0,0)=0$ by
continuity of $f$. Consequently, $\mathscr K[F](0)=\{0\}$.}

 {For stability,  {choose $\varepsilon>0$ such} that
$\overline{B}_\varepsilon(0)\subseteq\Rap$,  {take}
$c:=\min_{\|x\|=\varepsilon}\Vap(x)>0$ and
$\Omega:=\{x\in B_\varepsilon\mid\Vap(x)<c\}$, and choose $\delta>0$ with
$B_\delta\subseteq\Omega$, which is possible since $\Vap$ is continuous and
$\Vap(0)=0$. By Lemma~\ref{lem: recFeasInfMPC}, every maximal solution  {$\xx = \phi(\cdot, \uuinf, t_0, x_0)$ with}
$x_0\in B_\delta$ is complete, $\Rap$-valued, and locally absolutely
continuous, so \eqref{eq:V-decrease-inherited} yields
$\tfrac{\dd}{\dd t}\Vap( {\xx}(t))\le-W( {\xx}(t))\le0$ for almost every  {$t \geq t_0$} and
hence, $\Vap( {\xx}(t))\le\Vap(x_0)<c$ for all $t\ge t_0$. Such a solution  {therefore satisfies} $\| {\xx}(t)\|<\varepsilon$ for all
$t\ge t_0$.}

 {For convergence, let $x_0\in\Rap$ be arbitrary. Integrating
{the same time-differential form of \eqref{eq:V-decrease-inherited}} gives
$\int_{t_0}^{\infty}W(\phi(t))\,\dd t\le\Vap(x_0)<\infty$  {since $\Vap$ is nonnegative}. Every solution is
$M$-Lipschitz in $t$ and $W$ is uniformly continuous on the compact set $\Rap$,
so $t\mapsto W( {\xx}(t))$ is uniformly continuous and Barbalat's  {L}emma
\cite{khalil2002} gives $W( {\xx}(t))\to0$. Since $W$ is positive definite and
$\Rap$ is compact, this implies $ {\xx}(t)\to0$ as $t\to\infty$.}
\end{proof}

\subsection{QP-based Feedback Law}
 {In this subsection, we provide a practical formulation how the $\partial$MPC can be implemented, bridging the gap to standard  {QP-based} CLF-CBF approaches  {\cite{ames_control_2017, schneeberger_advanced_2026}} and showing how the efficient QP can be used in an MPC-like scheme.}
From here on, we assume that the stage cost $L$ is quadratic (in $u$) and that the system \eqref{eq:system} is control-affine.
 {Recall that the (approximation of) the feasible set is forward invariant \cite{chen1998} and that a sufficient condition for forward invariance is given by  {a} CBF. Likewise, the optimal value function or its approximation is a CLF \cite[Section 5.3]{gruenePannek2017}.}
Let $(\hat V, \hat h)$ be a  {any} compatible CLF/CBF pair that satisfies $\dot{\Vap}(x,\kappa(x))\leq-L(x,\kappa(x))$ and  {\eqref{eq:terminal-classK}} 
for some state feedback $\kappa$  {and extended class $\mathcal K$ function $\alpha$; that is}, the stabilizing terminal condition of Definition~\ref{def:terminal} are met  {by $\hat V$ and $\Rap :=\{x\in\mathbb{R}^n~|~\hat h(x) \leq 0\}\subseteq\CS$}.
The infinitesimal MPC problem for $\uuinfhat: \Rap\to\Uset$ is then given by the QP
\begin{subequations}
    \label{eq: infQP}
\begin{align}
    \hat\kappa_\partial: x \mapsto {} &\operatorname*{arg\,min}_{u \in \mathcal{U}} \; L(x, u) +  \ip{\nabla \Vap(x), f(x,u)}  \\
    &\text{subject to $\ip{\nabla \hat h(x), f(x,u)} \leq \alpha( -\hat h{(x)} )$}
\end{align}
\end{subequations}
for any given, sampled state $x \in \Rap$. The solution of ~\eqref{eq: infQP} preserves the stability and safety guarantees of the infinitesimal MPC problem~\eqref{eq:mpc-infinitesimal}, given by Theorem~\ref{thm:main}.  {In simple terms, a computationally lightweight scheme  {is} established which fuses the efficient computations of convex optimization with the theoretical guarantees of quasi-infinite horizon MPC.}

 { {Unlike} standard nonlinear discrete-time MPC with horizon $ {N} \geq 2$, where $T = N\delta$ with $N \in \mathbb{N}$ and $\delta >0$ is the sampling period (cf. \cite[p. 45]{gruenePannek2017}), which generally yields a nonconvex nonlinear optimization problem (NNOP), the proposed continuous-time  {$\partial$}MPC law is a convex QP, beneficial for online optimization. Table~\ref{tab:comparsionNMPCInfMPC} compares problem sizes of discrete-time NMPC, standard  {QP-based} CBF-CLF  {approaches}, and the proposed $\partial$MPC. The NMPC has extra equality constraints from the initial value problem; for inequalities, we assume simple upper and lower bounds, a terminal set constraint, and potential path constraints. Clearly, the QP {-based} approaches are drastically smaller -- even for $N = 1$, the NNOP is larger and still nonconvex. Unlike NNOP, for convex optimization problems it is possible to provide a polynomial upper bound on the number of iterations depending on the problem size (see \cite[p. 46]{malyuta2022} or \cite[Section 1.3.1]{boyd2004}). A small-sized problem results in less arithmetic operations and such an upper bound is beneficial for real-time applications. Note that the goal of this paper is not to establish such an upper bound  {but} an efficient scheme for which such an upper bound could be provided.}

\begin{table}[h!]
    \centering
    \caption{ {Size comparison regarding number of decision variables ($n_\mathrm{dec}$), equality constraints ($n_\mathrm{eq}$) and inequality constraint ($n_\mathrm{ineq}$) with potential path constraints ($n_\mathrm{path}$) between the proposed $\partial$MPC,  {QP-based} CBF-CLF  {approaches}, and discrete-time NMPC using a multiple-shooting formulation.} }
    \label{tab:comparsionNMPCInfMPC}
    \begin{tabular}{lcc c}
        \hline
    \hline
     &  NMPC & CBF-CLF & $\partial$MPC \\
     \hline
     $n_\mathrm{dec}$ & $(N+1)n +Nm $ & $ m + 1$ & $m$\\
     $n_\mathrm{eq}$ & $(N+1)n$ & 0  & 0 \\
     $n_\mathrm{ineq}$ & $2n_\mathrm{dec} + (N+1)n_\mathrm{path} + 1$ & $2n_\mathrm{dec} + 2$ & $2n_\mathrm{dec} + 1$ \\
         \hline
         \hline
    \end{tabular}
\end{table}

 {Both the standard  {QP-based} CBF-CLF  {approach} and the proposed $\partial$MPC  {scheme} result in small-sized convex QP, with the CBF-CLF approach slightly larger. To ensure feasibility, both employ a CBF as constraint. However, unlike 
CBF-CLF controllers  {such as \cite{ames_control_2017, schneeberger_advanced_2026}}, the proposed formulation allows to incorporate  {a} state-dependent stage cost, which allows to tune the controller performance. Whereas the $\partial$MPC ensures stability using a terminal penalty, the CBF-CLF  {approach} needs a CLF as constraint. To ensure feasibility, the CBF-CLF  {approach} needs an additional slack variable as decision variable that may relax the CLF constraint to (temporally) sacrifice stability for safety. Note that the CLF in the $\partial$MPC  {scheme} is not a standard CLF, but is determined such that the augmented cost function is an upper bound on the infinite horizon cost.}
 {\begin{remark}
MPC schemes without terminal constraints \cite{limon2006,gruene2012} or terminal conditions \cite[Chapter 6]{gruenePannek2017} exist with similar closed-loop properties. Note that the latter requires a sufficiently large prediction horizon and may become infeasible even if initially feasible~\cite[Section 7.4]{gruenePannek2017}.
\end{remark}}
 {To ensure that Theorem~\ref{thm:main} holds, the approximations $\Vap$ and $\hat h$ must satisfy the stabilizing terminal condition. Hand-crafted CLF/CBF pairs are cumbersome to find. Formal methods like Hamilton-Jacobi (HJ) reachability or SOS programming are suitable alternatives. HJ handles broad nonlinear systems but yields discontinuous feedback due to nonsmooth viscosity solutions and grid-dependent accuracy~\cite{gong_constructing_2023}. In contrast, we propose an SOS program to synthesize CBFs and CLFs via dissipation inequalities, avoiding spatial discretization. This yields smooth, twice-differentiable polynomial solutions needed for our scheme, and it offers a practical trade-off between complexity and conservatism~\cite{tan2006}}.


\section{Polynomial Synthesis}
\label{sec: SOSApproxAndQP}
In this section, we leverage sum-of-squares programming to compute polynomial approximations $\hat V$ and $\hat {\mathcal R}$ subject to the stabilizing terminal condition of {Definition~\ref{def:terminal}}. To that extent, we assume that the system dynamics $f$ in \eqref{eq:system} and the stage cost $L$ in \eqref{eq:mpc-problem} are polynomial functions, that the state constraint set (or a tight inner approximation) is the semialgebraic set
\begin{align*}
    \CS = \{x \in \mathbb R^n ~|~ g(x) \leq 0 \}
\end{align*}
for some $g \in \mathbb R[x]$, and that the input constraint set is the polyhedral set
\begin{align*}
    \Uset = \{ u \in \mathbb R^m ~|~ H_\Uset u \leq \mathbf 1_p \}
\end{align*}
with $H_\Uset \in \mathbb R^{p \times m}$, $p \in \mathbb N$, and $\mathbf 1_p$ denotes the vector $(1,\ldots,1) \in \mathbb R^p$. SOS methods have previously been used to approximate {reach-avoid set}s and optimal value functions (see, e.g., \cite{Cunis2021aut} and references herein); however, these works required polynomials that were parametrized in time. In our experience, these parameterizations lead to overly conservative results. Instead, we propose a nonconvex SOS optimization problem that directly approximates $V$ and $\mathcal R$ subject to the stabilizing terminal condition.

\subsection{Sum-of-Squares Programming}
Denote the set of polynomials in $x$ with real coefficients up to degree $d$ by $\mathbb R_d[x]$. 
A polynomial $p \in \mathbb R_{2d}[x]$ is a {\em sum-of-squares polynomial} ($p \in \Sigma_{2d}[x] \subset \mathbb{R}_{2d}[x]$) if and only if there exist $m \in \mathbb N$ and $p_1, \ldots, p_m \in \mathbb R_d[x]$ such that $p = \sum_{i=1}^m (p_i)^2$. If $p \in  \Sigma_{2d}[x]$, then $p(x) \geq 0 $ for all $x \in \mathbb{R}^n$.
    
A nonconvex SOS problem takes the form~\cite{Cunis2025acc}
\begin{align}
    \label{eq:sos-nonlinear}
    \min_{\xi \in \mathbb R_{2d}[x]^n} f(\xi) \quad \text{s.t. $\xi \in \Sigma_{2d}[x]^n$ and $g(\xi) \in \Sigma_{2d'}[x]^m$}
\end{align}
where $f: \mathbb R_{2d}[x]^n \to \mathbb R$ and $g: \mathbb R_{2d}[x]^n \to \mathbb R_{2d'}[x]^m$ are differentiable functionals.

Unlike convex SOS program{s}, which can be transcribed into and solved as semidefinite programs~(SDP)~\cite{parrilo2003}, nonconvex problems need iterative schemes such as coordinate-descent~\cite{chakrabortyEtAl2011} or sequential SOS~\cite{Cunis2023acc}. Our recently introduced MATLAB-based software suite Ca$\Sigma$oS~\cite{Cunis2025acc} is specifically designed to solve such problems, whereas other toolboxes allow to manually implement iterative schemes using convex SOS programs.

\subsection{Sufficient Condition Estimation}
\label{subsec: SuffConEst}
As a first step in implementation, we are interested in synthesizing a polynomial $\hat h \in \mathbb R_{2d}[x]$ and a positive definite polynomial $\hat V \in \mathbb R_{2d}[x]$ that satisfy the stabilizing terminal condition with polynomial feedback $\hat\kappa \in \mathbb R[x]^m$; thus, $(\hat V, \hat h)$ are a compatible polynomial CLF/CBF pair  {that fulfills the conditions of Definition~\ref{def:terminal}  {with respect to the stage cost $L$ and input constraints $\Uset$}}.  
 {Unlike the standard CBF-CLF synthesis, we propose a nonconvex optimization framework to approximate the infinite-horizon controllability set and optimal value function of the MPC problem \eqref{eq:mpc-problem}.}
Our approach is similar to~\cite{schneeberger_sos_2023} where CLF and CBFs both have the so-called \textit{control-sharing property} to be compatible.
In the following, we provide a step-by-step \textit{recipe} of the synthesis problem which can be ultimately solved by SOS programming.
Note that the state dependence is omitted if clear from context.
\subsubsection*{Stabilizing Terminal Condition}
We enforce the set-containment conditions
\begin{subequations}
    \label{eq: innerApproxCBF}
\begin{align}
    &\{x \in\mathbb R^n \mid \hat h(x) \leq \beta\} \subseteq \{x \in \mathbb{R}^n \mid g(x) \leq 0\} \\
    &\{x \in \mathbb R^n \mid \hat h(x) \leq \beta\} \subseteq \{x \in \mathbb{R}^n \mid  \dot{\hat{h}}(x) \leq \alpha (\beta -\hat h(x) )\}  \label{eq: dissipatIonInv}
\end{align}
\end{subequations}
with $\dot{\hat{h}}(x) = \ip{\nabla \hat h(x), f(x,\hat\kappa(x))}$,
where 
$\alpha( \cdot )$ is an extended class $\mathcal{K}$ function
to ensure forward-invariance (cf.~\cite[Section~III]{ames_control_2019}) and constraint satisfaction of the feasible set approximate,  
\begin{align*}
    \mathcal{\hat{R}}:= \{x \in \mathbb R^n \mid \hat h(x) \leq \beta \} \subseteq \CS
\end{align*}
for some level $\beta > 0$.

We additionally impose~\cite{jarvisEtAl2003}
\begin{align}
    &\{x \in \mathbb R^n \mid \hat h(x) \leq \beta\} \subseteq \{x \in \mathbb{R}^n \mid  H_\mathcal{U} \hat\kappa(x) \leq \mathbf 1_p\} \label{eq: sosContCon2}
\end{align}
to ensure that the synthesized control law is admissible.


To ensure that $\hat V$ satisfies condition \eqref{eq:terminal-cost} over $\hat {\mathcal R}$, we make use of the continuous-time dissipation inequality 
\begin{align*}
    \tau(x) := \langle \nabla \hat V(x), f(x,\hat\kappa(x)) \rangle + L(x,\kappa(x)) \leq 0 
\end{align*}
and the set-containment condition
\begin{align}
     &\{x \in \mathbb R^n \mid \hat h \leq \beta\} \subseteq \{x \in \mathbb{R}^n \mid \tau(x)\leq 0\} \label{eq: penaltyDissip}
\end{align}
for $\hat V$, $\hat h$, and $\hat \kappa$.

\subsubsection*{Synthesis Problem}
To reduce conservatism in terms of size of the approximate feasible set $\hat {\mathcal R}$, we want to synthesize an invariant set which is as close as possible to the boundary of the constraint set. One way to achieve this is to minimize the squared $l_2$-norm $\|\cdot\|_{\mathbb R[x]}$ on the space of polynomials. 
Potentially, one is also interested in a good approximation of the optimal value function. To balance both goals, we propose the multiobjective cost function
\begin{align}
    J(\hat V, \hat h, \hat\kappa) = \lambda_1 ||g-\hat h||_{\mathbb R[x]}^2 + \lambda_2|| \ip{\nabla \hat V, f} - L||_{\mathbb R[x]}^2 \label{eq: multiObjCost}
\end{align}
where $\lambda_1 + \lambda_2 = 1$ are user-specified weights.
The high-level optimization problem reads
\begin{subequations}
\begin{align*}
    &\text{minimize} \; J(\hat V, \hat h, \hat\kappa) \\
    &\text{subject to \eqref{eq: innerApproxCBF}--\eqref{eq: penaltyDissip}}
\end{align*}
\end{subequations}

\subsubsection*{SOS Program}
To transcribe the high-level optimization problem into a SOS problem we make use of the so-called generalized $\mathcal{S}$-procedure, a relaxed version of the Positivstellensaz~\cite{wang2023}, as outlined next.
\begin{theorem}\cite[Lemma~2.1]{tan2006}
    Given $p_0, p_1, \ldots, p_N \in \mathbb{R}[x]$, if there exist $s_1, \ldots, s_N \in \Sigma[x]$ such that $p_0- \sum_{k=1}^N s_k p_k \in \Sigma[x]$, then
	\begin{align}
		\bigcap_{k=1}^N \{x \in \mathbb{R}^n \mid p_k(x) \geq 0\} \subseteq \{x \in \mathbb{R}^n \mid p_0(x) \geq 0\}. \nonumber
	\end{align}
    \label{Theo: genSprocedure}
\end{theorem}
Applying the generalized $\mathcal{S}$-procedure to the constraints \eqref{eq: innerApproxCBF}--\eqref{eq: penaltyDissip} yields the following SOS problem
\begin{subequations}
    \label{eq: SOSprob1}
\begin{alignat}{2}
    &\min_{\substack{s_1, s_2, s_4 \in\Sigma[x], s_3 \in\Sigma[x]^p \\  \hat h , \hat V \in \mathbb{R}[x], \hat\kappa \in \mathbb{R}[x]^m} } \; J(\hat V, \hat h, \hat\kappa)   \\
    &\text{subject to }  \hat V - \varepsilon (x^\top x)                &  & \in\Sigma[x]^{}  \label{eq: strictPos}\\
    &\hphantom{\text{subject to }}   s_1 (\hat h-\beta) - g                                     &  & \in \Sigma[x]^{} \\
    &\hphantom{\text{subject to }} s_2 (\hat h-\beta) -  \dot{\hat{h}} - \alpha ( \hat{h} - \beta )&  & \in \Sigma[x]^{} \label{eq: CBFSOS}\\
    &\hphantom{\text{subject to }}  s_3 (\hat h-\beta) - (H_\mathcal{U} \hat\kappa - \mathbf 1_p)                 & & \in \Sigma[x]^p \label{eq: lowerContSOS} \\
    &\hphantom{\text{subject to }}  s_4 (\hat h-\beta) - \tau                                    &  & \in\Sigma[x]^{} \label{eq: terminalPenSOS}
\end{alignat}
\end{subequations}
where~\eqref{eq: strictPos} ensures positive definiteness of $\hat V$ with a small constant $\varepsilon > 0$ and in~\eqref{eq: CBFSOS}, the extended class $\mathcal{K}$ function is of the form $\alpha(s) = a \, s$ with $a > 0$.
This problem is nonconvex and cannot be directly solved by an SDP. Instead, we make use of a sequential SOS approach~\cite{Cunis2023acc} implemented in Ca$\Sigma$oS~\cite{Cunis2025acc}. The latter one allows to directly optimize the nonlinear problem instead of solving a sequence of alternating subproblems. Once the SOS synthesis step is done, we store $\hat V$ and $\hat h - \beta$ for online usage.

\section{Numerical Results}
\label{sec: NumericalResults}
In this section we evaluate the proposed approach in numerical examples from nonlinear spacecraft control. Compared to Earth-based applications, computing resources are severely limited by space-qualified hardware~\cite{eren_model_2017}. Thus, approaches that reduce online computational effort are paramount.  We consider two case studies, a large-angle slew maneuver and a collision avoidance scenario, of a satellite similar to the Hubble telescope~\cite{nurre1995}.
For simplicity, the dynamics in the simulations and the optimal control problem are identical.

For the precomputation step as described in Section~\ref{subsec: SuffConEst}, we make use of {\rm Ca$\Sigma$oS}~\cite{Cunis2025acc} v1.0.0-rc with {\rm MOSEK}~\cite{andersen_mosek_2000} v11.0.28 as the underlying SDP solver. For online optimizations, we make use of {\rm CasADi}~\cite{andersson_casadi_2019} v3.6.7 to setup the problems. All solvers are set up using their default options but we turn off display outputs. If not stated otherwise, computation time (wall-time) is provided by the CasADi solver statistics.
The results in this paper were computed on a  computer with Windows 10 and MATLAB 2023b running on an AMD Ryzen 9 5950X 16-Core Processor with 3.40 GHz and 128 GB RAM. Source code for the synthesis  and the simulations along with details of the polynomial optimization problem \eqref{eq: SOSprob1} is provided in the supplementary material~\cite{darusSupplemenataryMaterial}. 

\subsubsection*{Dynamic and Constraints}
To describe the satellite attitude, we make use of Modified Rodrigues Parameters (MRP) denoted by $\sigma \in \mathbb R^3$, which allows for potentially large angle slews, although it has a singularity at $2 \pi$.  We consider $\sigma^T \sigma \leq 1$, which guarantees that all rotations are shortest-path~\cite{junkins_analytical_2009}, at the cost of restricting to only a subset of potential singular free re-orientation maneuver. 

The state vector reads $x = (\omega, \sigma)$ where   $\omega \in \mathbb R^3$ are the angular rates.
The system dynamics  $f: \mathbb R^6 \times \mathbb R^3 \rightarrow \mathbb R^6$ in continuous-time read
\begin{subequations}
    \label{eq:spacecraft-dynamics}
\begin{align}
\dot{x} = f(x,u) = \begin{bmatrix}  -J^{-1} \tilde{\omega} J \omega + J^{-1} u  \\  \frac{1}{4} B(\sigma) \omega \end{bmatrix} 
\end{align}
with $B(\sigma) = (1-\sigma^T \sigma) I_{3\times 3} + 2 \tilde \sigma + 2 \sigma \sigma^T$,
where $J \in \mathbb R^{3 \times 3}$ is the inertia tensor, $u \in \mathbb R^3$ are the control torques, and $\tilde{\omega}$ and $\tilde{\sigma}$ are the cross-product matrix of the angular rates and MRP respectively. We consider the inertia tensor of the Hubble telescope that reads $J = \text{diag}(31046, 77217, 78754)$ \unit{\kilogram\meter{^2}}~\cite{nurre1995}. The following state and control constraints are considered
\begin{align}
    \omega &\in [-0.5, 0.5] \times  [-0.2, 0.2] \times [-0.2, 0.2]~\si{ \unit{\deg\per\second}} \label{eq:rateCon}\\
    \sigma^T \sigma &\leq 1 \label{eq: MRPConst}\\
    \| u \|_\infty &\leq \SI{1.2}{\newton\meter}
\end{align}
\end{subequations}

\subsection{Comparative Study}
\label{subsec: consSatSimple}
We compare our $\partial$MPC with polynomial approximations of optimal value function and feasible set to the discrete-time  nonlinear MPC (NMPC) formulation \cite[p.~100]{gruenePannek2017}  {with a prediction horizon $N \geq 2$}, an RTI scheme,\footnote{We make use of CasADi's {\tt sqpmethod}, which we limit to one iteration (cf.~\cite{diehl_real-time_2005,frey2024}) and provide the previous solution as initial guess.} the polynomial control law  $\hat\kappa(x)$ from the synthesis step in \eqref{eq: SOSprob1}, a standard QP-based CBF-CLF approach   {(see \cite{darusSupplemenataryMaterial})}, for which we precompute a compatible CLF/CBF pair offline using an adapted version of~\eqref{eq: SOSprob1} 
 {and a discrete-time NMPC with prediction horizon $N = 1$ (horizon-one), which also represents a computational lower bound for tracking MPC with artificial reference, e.g., as in \cite{kohler2020b}. For simplicity, the horizon-one scheme uses the same ingredients as the proposed $\partial$MPC but in discrete-time (see \cite{darusSupplemenataryMaterial}).}
{All three $\partial$MPC, RTI and CBF-CLF require a QP to be solved during simulation, for which we used {\rm qrqp}~\cite{ANDERSSON2018331}. For the NMPC  {and horizon-one NMPC}, we use {\rm Ipopt}~\cite{wachter_implementation_2006} to solve the discrete-time OCP. We also tried to use the more efficient solvers {\rm alpaqa}~\cite{pas2022a} and {\rm fatrop}~\cite{vanroye2023fatrop}; unfortunately, these solvers either failed due to infeasibility or took significantly longer than {\rm Ipopt} to solve the problem and are hence omitted.

\subsubsection*{Problem Description}
We assume a rest-to-rest profile, meaning the initial and final angular rates are zero. To compare the five approaches, we consider three single-axis large angle slew maneuvers with rotations of $\chi \in {} \{75^\circ,90^\circ, 110^\circ\}$. 
We use a uniform discretization with $\Delta t = T/100$ for all three cases. We heuristically determine the prediction horizon $T$  for NMPC and RTI such that the first step is feasible. We find $T \in {}\{\SI{200}{\second},\SI{300}{\second}, \SI{400}{\second}\}$ to be feasible, which can be primarily traced back to the large entries in the inertia tensor.
In the MPC formulations, we use the stage cost $L(x,u) = x^T Q x + u^T R u$ with matrices $Q = I_{6\times6}$ and $R = I_{3\times3}$, where $I$ denotes the identity matrix. The QP-based CBF-CLF approach uses the same $R$ matrix.

\subsubsection*{Test Case Preparation}
For our $\partial$MPC approach,  we solve \eqref{eq: SOSprob1} to precompute the polynomial approximations $\hat V$ and $\hat h$ where the above described constraints are encoded in an inner-approximation of the constraint set   {(see \cite{darusSupplemenataryMaterial})}and the weights $\lambda_1 = 1$ and  $\lambda_2 = 0$ are used in \eqref{eq: multiObjCost}. We make use of a quadratic polynomial for the CBF and quartic for the CLF. We heuristically determined $a = 0.0001$ and fixed $\beta = 0.9$ for the synthesis step.
For the NMPC and RTI schemes, we determine a terminal penalty for the linearized system $F: x \mapsto x^\top S x$, where $S \in \mathbb R^{6 \times 6}$ is the solution of the Riccati equation for $Q$ and $R$. We find the maximum stable level set $\gamma  = 60.5$ for the nonlinear constrained system using a SOS program via bisection (see \cite{darusSupplemenataryMaterial} for more details) and define $\TS = \{ x \in \mathbb R^6 \, | \, F(x) \leq \gamma \}$. 
For  {settling time} analysis and visualization, we transcribe the MRP $\sigma$ into Euler angles
$\Phi = (\phi, \theta, \psi)$. We consider the dynamics as converged if
$ \lVert \omega(t) \rVert_\infty \leq 0.001~\si{\deg\per\second}$, $\lVert \Phi(t) \rVert_\infty \leq 0.3~\si{\deg}$, and $\lVert u(t) \rVert_\infty \leq 0.001~\si{\newton\meter}$. The maximum simulation time is set to \SI{5000}{\second} with a fixed step size of $\SI{0.1}{\second}$ ($\SI{10}{Hz}$).

\subsubsection*{Results}
We compare the  {six} approaches in four categories aggregated over the three initial conditions in Table~\ref{tab:comparsionSatSimple}. We evaluate mean  {settling time} over three samples, integral stage cost,\footnote{We evaluate $\int\! L(\cdot) \mathrm d t$ along the simulated, closed-loop state and input trajectories.} mean and worst-case computation times, and precomputation time for each approach, marking the best in \textbf{bold}. Detailed results are provided in the supplementary material~\cite{darusSupplemenataryMaterial}. Fig.~\ref{fig:satSimpleRotXaxis} shows roll-rate, roll angle, and control input for all approaches during the simulation for the $\chi = {}110^\circ$ maneuver. All approaches are feasible but exhibit different convergence behaviors. Performance highly depends on cost-function weights and other parameters and can be improved by adequate tuning.

The $\partial$MPC approach performs best in  {settling time}, integral stage cost, and computation time. The QP-based approaches outperformed NMPC schemes in online computation time, which is expected. The polynomial control law from the synthesis step requires no online optimization and therefore, no computation time is reported.  {
The horizon-one NMPC performs mildly worse in settling time with the same integral stage cost, but its computation times are several magnitudes larger than those of the QP-based approaches.}

The computation time for NMPC exceeds the sampling frequency, making real-time application impossible. The RTI scheme is faster than NMPC, coming closer to real-time execution despite some worst-case computations exceeding the sampling frequency.
Fast online computation times of the QP-based approaches,  {the horizon-one scheme} and the polynomial control law $\hat\kappa(\cdot)$ stem from offline precomputation, which  took slightly more than a minute. In contrast, MPC and RTI have shorter precomputation times.

MPC and RTI had similar but worse  {settling times} than the proposed approach. The QP-based CBF-CLF and polynomial law performed worst, with the CBF-CLF approach slightly better. Performance might improve with different comparison functions or hand-crafted CLF/CBF.
Furthermore,
the $\partial$MPC approach's integral cost is compatible with MPC and RTI. CBF-CLF  ranks third, slightly worse despite longer  {settling times}, while considering only control in its cost. The polynomial law, though less competitive, ensures stability and safety without online optimization.

\begin{figure}[h!]
    \centering
        \setlength{\figH}{4.1cm}
    \setlength{\figW}{7cm}
     \input{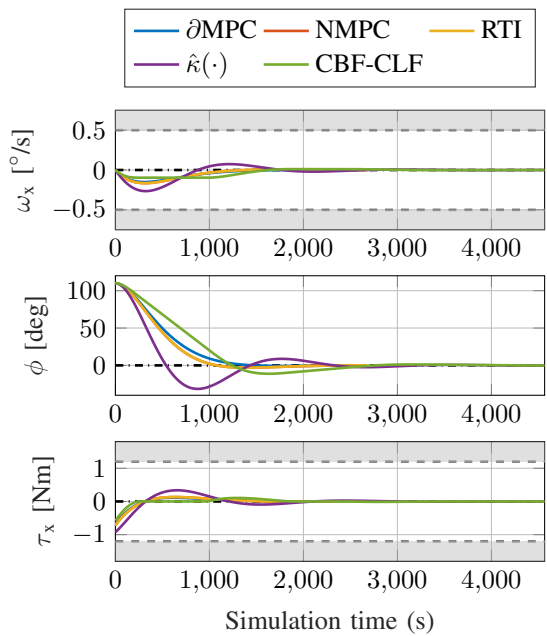}
    \caption{Comparison of $\partial$MPC,  NMPC, RTI scheme, polynomial control law, QP-based CBF-CLF,  {and horizon-one MPC} for the $\chi = \ang{110}$ single-axis re-orientation maneuver.}
    \label{fig:satSimpleRotXaxis}
\end{figure}

\begin{table}[h!]
    \centering
    \caption{Aggregated comparison of the different controllers in terms of precomputation time $t_\text{pre}$, mean  {settling} time $t_\text{set}$, mean (worst-case) computation time $t_\text{comp}$, and mean integral stage cost for three different initial conditions.}
    \label{tab:comparsionSatSimple}
    \begin{tabular}{lccrlc}
        \hline
    \hline
    Method & $t_\text{pre}$ [s]& $t_\text{set}$ [s] & \multicolumn{2}{c}{$t_\text{comp}$ [ms]} & $\int L(\cdot) \mathrm dt$ \\
    \hline
     $\partial$MPC      &  72.3 & $\mathbf{2275.8}$  &  $\mathbf{0.003}$ & $\mathbf{(0.02)}$   & $\mathbf{88.8}$\\
     NMPC ({\rm Ipopt})  &   $\mathbf{16.4}$    & 2586.5  & 1009.4 & (11842) &  89.8\\
     RTI           &   $\mathbf{16.4}$    &  2586.6 & 51.2 & (585.6) & 89.8\\
     $\hat\kappa(\cdot)$     & 72.3  &  4297.2 & {---} & ({---})    & 145.7\\
     CBF-CLF    &  70.2 &  4163.8 & $\mathbf{0.003}$ & ${(0.1)}$     & 96.02\\
     { Horizon-one}    &   {72.3} &   {2276.5} &  {$3.8$} &  {($11.8$)}     &  {\textbf{88.8}}\\
         \hline
         \hline
    \end{tabular}
\end{table}

\subsection{Performance Test}
We now evaluate our $\partial$MPC approach in a challenging task for spacecraft control, where the reorientation is subject to mission-critical state constraints~\cite[Section~3.4]{malyutaAdvancesTrajectoryOptimization2021a}. As we demonstrate, our precomputed polynomial approximations of optimal value function and feasible set cover large portions of the state space. A similar coverage would require long prediction horizons in NMPC or RTI schemes, compromising the online computation times.

\subsubsection*{Problem Description}
A difficult task for scientific observation spacecraft is to perform large angle three-axis slew maneuvers while ensuring safety (non-exposure) of sensitive on-board instruments such as telescopes \cite[Section~3.4]{malyutaAdvancesTrajectoryOptimization2021a}. Such strict nonexposure (attitude pointing) constraints are frequently described by conic state constraints of the form (cf.~\cite{diazramosKinematicSteeringLaw2018})
\begin{align}
    \langle {_I}n, T_{IB}({_B}b) \rangle - \cos(\delta_\text{min}) < 0
    \label{eq: keepOutCone}
\end{align}
where ${_I}n \in \mathbb R^3$ is a unit vector in the inertial frame describing the direction of the cone to be avoided,  $T_{IB}: \mathbb R^3 \rightarrow \mathbb R^3$ is the linear transformation from body to inertial coordinates, ${_B}b \in \mathbb R^3$ is the body-fixed boresight vector of an instrument, and $\delta_\text{min} > 0$ describes the safety angle. The main challenges in this problem arise from the nonlinear, nonconvex spacecraft dynamics and nonconvex conic constraints \cite[Section~3.4]{malyutaAdvancesTrajectoryOptimization2021a}.

We add one keep-out cone constraint~\eqref{eq: keepOutCone} to the problem description given in \eqref{eq:spacecraft-dynamics}, with ${_I}n = (0, 1, 0)$, ${_B}b = (1, 0, 0)$ and $\delta_\text{min} = \SI{20}{\deg}$. 
Additionally, we relax \eqref{eq: MRPConst} to a larger level set of six (i.e., $\sigma^\top \sigma \leq 6)$. This increases the domain of attraction of the proposed approach as larger slew maneuvers are allowed while ensuring singularity free motions. This comes at the cost of non-shortest path, i.e., the controller might take unnecessarily long paths.  Again, we are interested in an rest-to-rest maneuver. 

\subsubsection*{Test Case Preparation}
For the keep-out cone constraint, relaxed MRP constraint, and rate constraints we compute an inner-approximation using the first approach in \cite{darusSupplemenataryMaterial}.  We use the same cost function weights as in the previous test case and precompute the polynomial approximations solving \eqref{eq: SOSprob1}. Here, we only aim to reduce the distance between the CBF and the inner-approximate of the constraint set for angular rates set to zero in~\eqref{eq: multiObjCost}. We make use of a quartic CBF and a quadratic CLF in this example. We set $\beta = 0.01$ and make use of a linear comparison function $\alpha$ with $a = 0.001$.

The precomputation of the inner approximation took about \SI{35}{\second} and the precomputation of the CLF and CBF about \SI{236}{\second} with a good initial-guess. We perform a Monte-Carlo simulation campaign for three-axis constrained satellite re-orientation (rest-to-rest) maneuvers. We draw 100 uniformly distributed samples of which the initial attitude lies in the feasible set, that is, $\hat h(x_\text{sample}) \leq 0$. 

\subsubsection*{Results}
\begin{figure}[t]
\centering
\begin{tikzpicture}
  \node[draw,inner sep=3pt] (L) {%
    \begin{tabular}{@{}lc@{\qquad}lc@{\qquad}lc@{\qquad}}
      \legendline{red}        & ${_I}e_x$   &
      \legendline{green}      & ${_I}e_y$  &
      \legendline{blue}       & ${_I}e_z$   \\
            \legendstar{blue}       &  origin &
      \legendstar{green}      & ${_I}b(0)$    &
      \legendline{green}      & ${_I}b(t)$     
    \end{tabular}%
  };
\end{tikzpicture}


    \hfill
\begin{subfigure}[b]{0.45\columnwidth}
  \includegraphics[width=\textwidth,clip,trim=50 120 30 120]{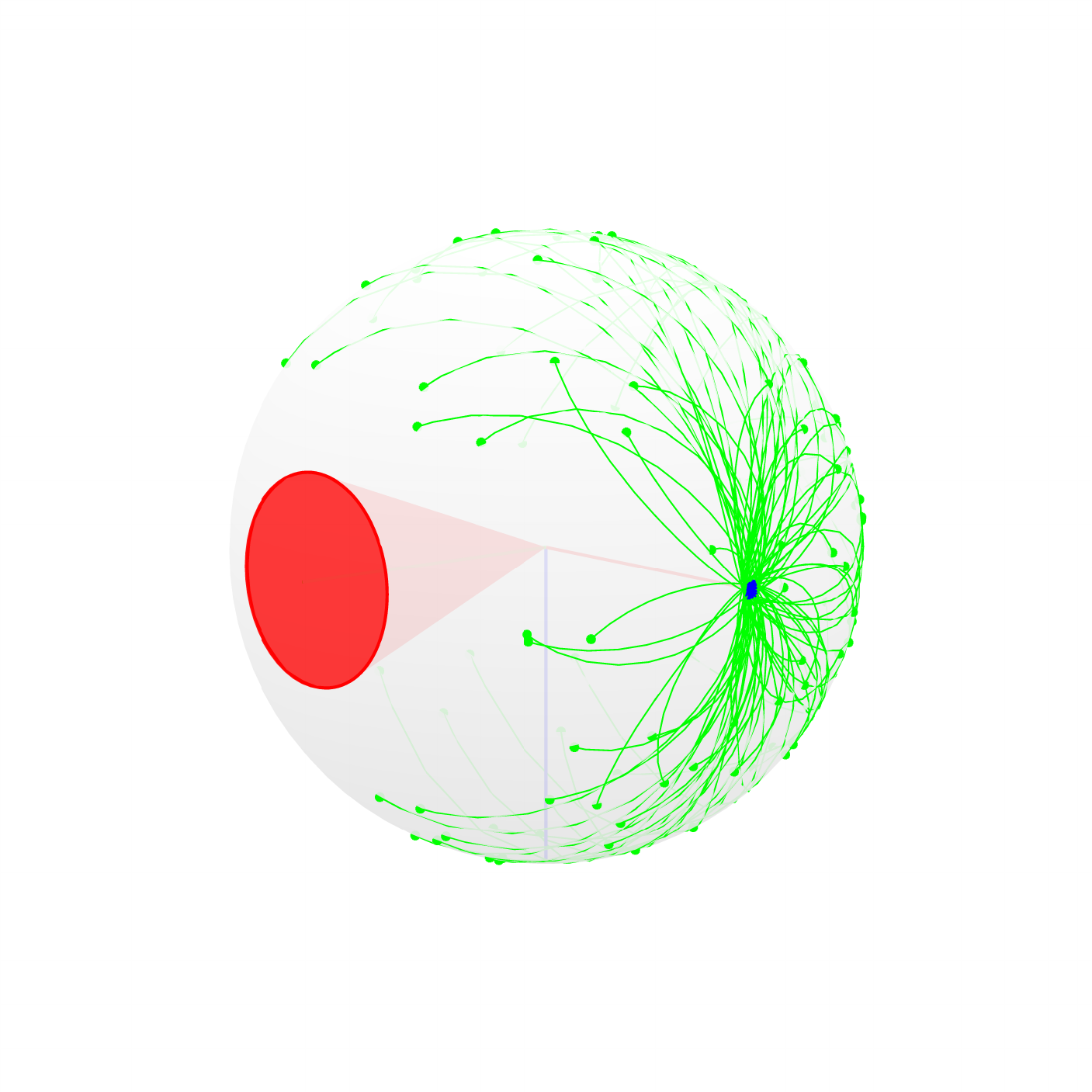}
\end{subfigure}
    \hfill
\begin{subfigure}[b]{0.45\columnwidth}
  \includegraphics[width=\textwidth,clip,trim=50 120 30 120]{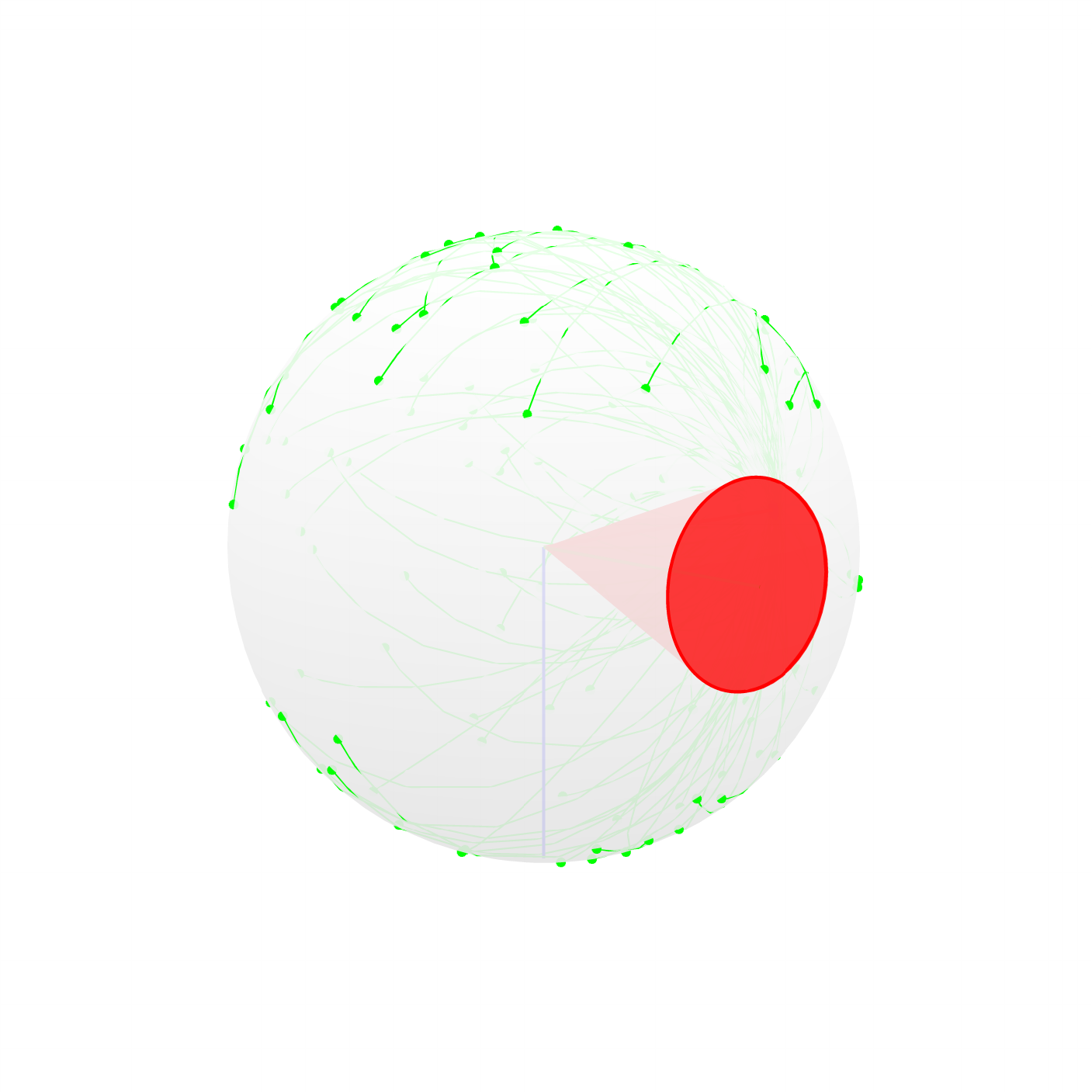}
\end{subfigure}
    \hfill

\caption{3D illustration of the three-axis constrained satellite re-orientation with one keep-out cone (red) on the unit sphere with inertial axes ${_I}e_{(\cdot)}$. Green stars mark initial conditions; green lines show the boresight trajectory ${_I}b(t)$ in inertial coordinates. All samples converge to the origin (blue star) while avoiding the keep-out cone.}
\label{fig:keepOutMC_unitSphere}
\end{figure}

The considered scenario is depicted in Fig.~\ref{fig:keepOutMC_unitSphere}. The plot shows the uniformly sampled initial conditions (green stars) and the corresponding trajectories (green lines).
All trajectories converge to the stable equilibrium (origin; blue star) while avoiding the {obstacle} (red cone). The feasible samples cover large portions of the unit sphere. This indicates a large feasible set, which is also the domain of attraction of the proposed approach. 
Fig.~\ref{fig:sat3D_MC_suffCon} depicts the evaluation of the approximated sufficient condition~\eqref{eq:terminal-cost} for asymptotic stability along the trajectories, which is fulfilled for all samples at any time. 
Similar, all samples that start in the feasible set, do not violate the keep-out cone constraint. Fig.~\ref{fig:sat3D_MC_barrier} depicts the polynomial CBF evaluated along the closed-loop trajectories. Some of the samples start fairly close to the boundary ($\hat h$ is almost zero), but never cross the zero line and are thus safe at all time. Clearly, none of the trajectories selects a path closely along the feasible set boundary. 
One the one hand, this can be traced back to samples that do start far from the boundary (e.g. the opposite site of the keep-out cone). On the other hand, trajectories that do not take a more direct \textit{shortest-path} can be traced back to relaxed MRP constraint (no guarantee for shortest-path) and to the  polynomial describing the feasible set. The latter can be improved by using different cost functions and higher-order polynomials during precomputation. The selection of $a$ in the comparison function may also affect the behavior.
The mean computation time was \SI{0.003}{\milli\second} and the worst case computation time~\SI{0.08}{\milli\second}. More plots are available in the supplementary material~\cite{darusSupplemenataryMaterial}.

 \begin{figure}[h!]
    \centering
    \setlength{\figH}{1.65cm}
    \setlength{\figW}{7cm}
    \input{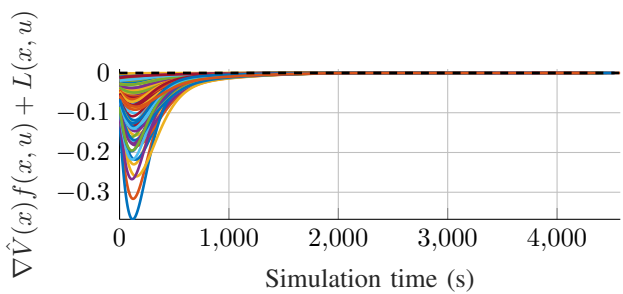}
    \caption{Sufficient condition~\eqref{eq:terminal-cost} evaluated along trajectories. Crossing the zero line corresponds to violation. The sufficient condition is satisfied at any time.}
    \label{fig:sat3D_MC_suffCon}
\end{figure}

\begin{figure}[htb!]
    \centering
    \setlength{\figH}{1.65cm}
\setlength{\figW}{7cm}
    \input{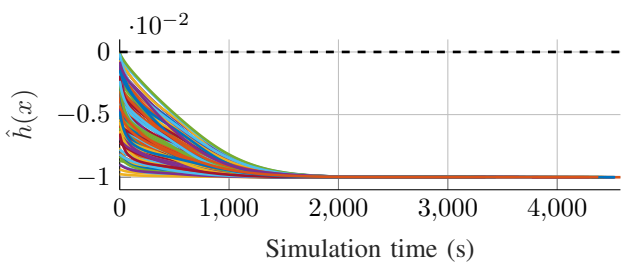}
    \caption{CBF evaluated along trajectories. Crossing the zero line corresponds to violation. All samples stay safe at any time.}
    \label{fig:sat3D_MC_barrier}
\end{figure}

\subsection{Discussion}
The $\partial$MPC approach guarantees asymptotic stability and constraint satisfaction by precomputed polynomial approximations. 
The solution quality (size of feasible set, control performance) highly depends on the polynomial degrees of the approximations. Thus, compared to classical NMPC, the offline tuning effort is higher because the polynomial structure (degree, monomials) must be selected in addition to the state cost and its weights. Selecting the extended class $\mathcal{K}$ functions affects the performance of both $\partial$MPC and {standard} CBF-CLF approaches and is an ongoing research activity (see, e.g.,~\cite{parwana2025}). The optimal value function is only approximated and hence one cannot expect it to always compete with or outperform, e.g.,  NMPC. The low computational effort of $\partial$MPC can be primarily traced back to the offline precomputation step. However, the pre-computed terminal ingredients are only valid for one equilibrium point. This is also the case for NMPC with stabilizing terminal condition.  Previous works on adapting a single terminal set \cite{Simon2014} or computing reference-dependence terminal sets \cite{cotorruelo_reference_2021} are applicable here, too. For the above satellite example, one could also consider a list of desired attitudes (set-points).  {The} $\partial$MPC  {is competitive to} established methods for nonlinear control  {with a significantly smaller computational footprint} and allows for tasks which would otherwise require prohibitively long prediction horizons.  {The $\partial$MPC is an advancement for MPC schemes with short prediction horizons and control-affine dynamics, i.e., a nonlinear MPC scheme formulated as a convex optimization problem.}

\section{Conclusion}
\label{sec: Conclusion}
 {This paper focuses on a computationally lightweight MPC scheme under state and input constraints. This approximate scheme fuses the computational efficiency of quadratic programs (QPs) that leverage control barrier functions (CBFs) and control Lyapunov function (CLF), with the theoretical backend of quasi-infinite horizon MPC. In numerical comparison to other control techniques for nonlinear systems including MPC, it is demonstrated that the proposed approach is competitive in performance while being significantly faster in computation time. The results indicate high online efficiency in terms of resource utilization that motivates future work on embedded hardware with limited processing capabilities, the latter which is common in aerospace applications.}

\bibliographystyle{IEEEtran}
\bibliography{references,MPC}
\begin{IEEEbiography}[{\includegraphics[width=1in,height=1.25in,clip,keepaspectratio]{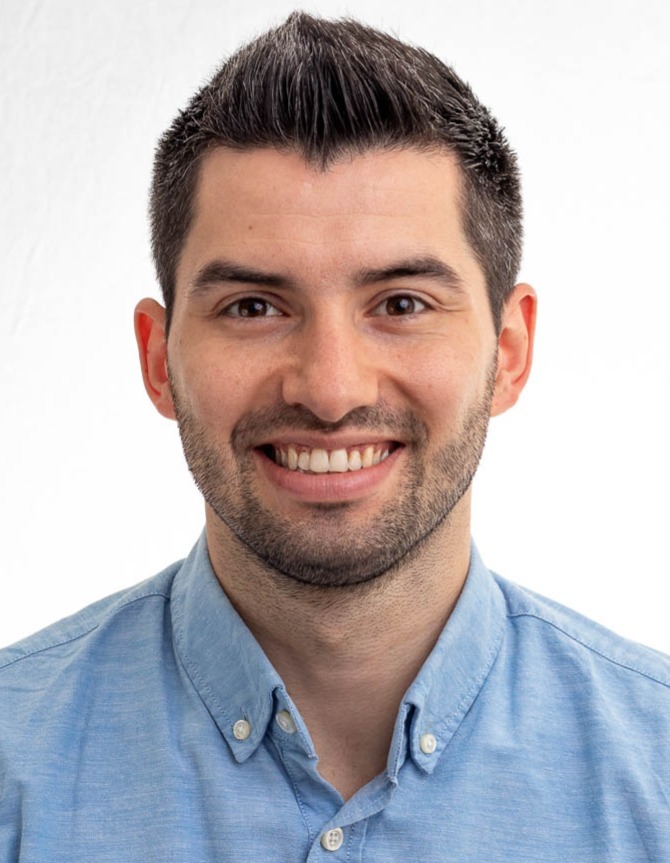}}]{Jan Olucak} received the M.Sc. degree
in aerospace engineering from the University of
Stuttgart, Stuttgart, Germany, in 2019.

Since 2019, he has been a Ph.D. Student
with the Institute of Flight Mechanics and Controls of University of Stuttgart, under supervision of Prof. Walter Fichter. His research interest lies in real-time optimal control methods.
\end{IEEEbiography}
\begin{IEEEbiography}[{\includegraphics[width=1in,height=1.25in,clip,keepaspectratio]{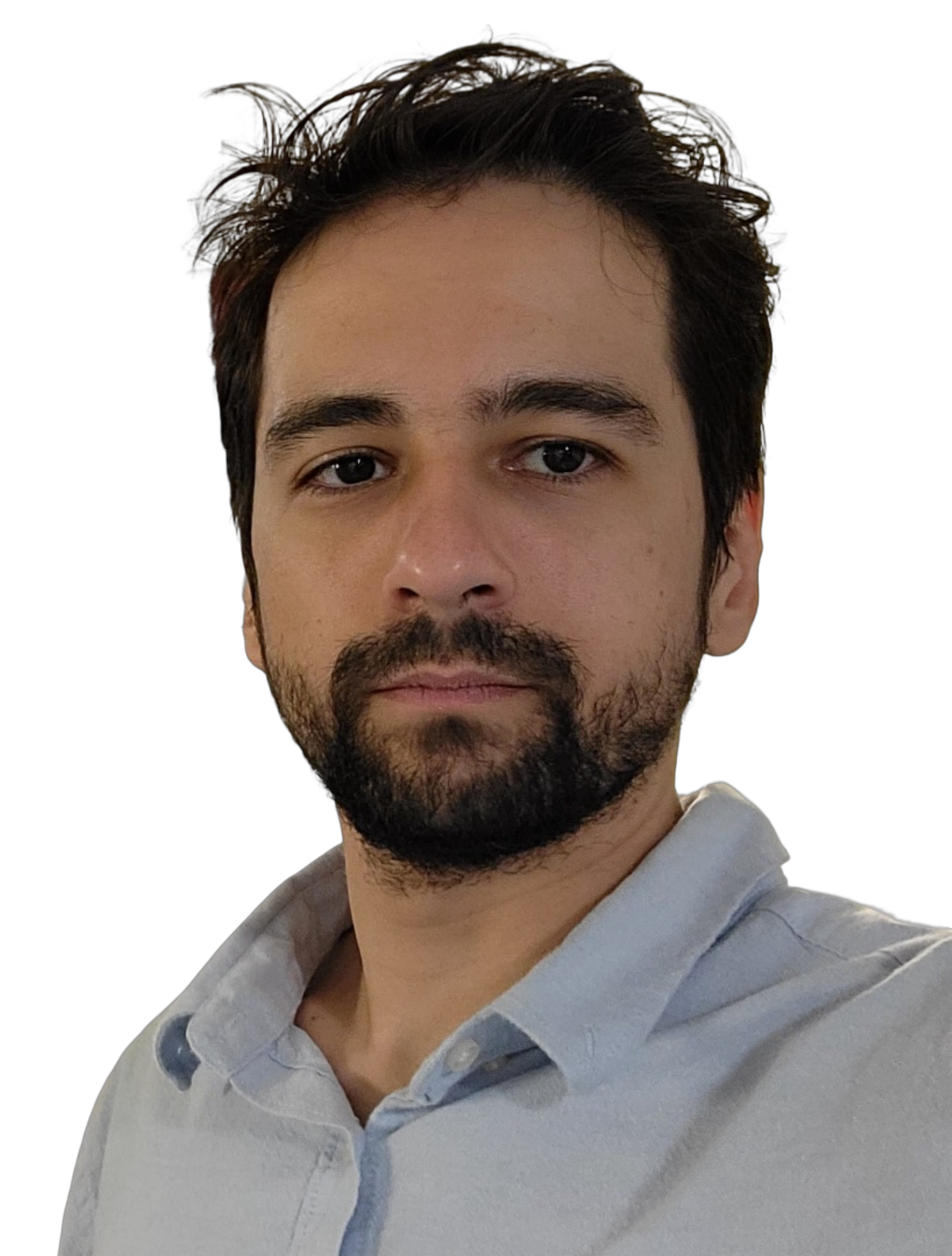}}]{Arthur Castello B. de Oliveira}
received his B.Sc. in electrical engineering with emphasis in control and M.Sc. in control engineering from the University of São Paulo in 2016 and 2019 respectively. He received his Ph.D. degree in electrical engineering from Northeastern University in 2024. He is currently a postdoctoral research associate at Northeastern University, Boston, MA.
His research interests include nonlinear systems and control, nonlinear optimization and optimal control, gradient systems, artificial intelligence techniques and machine learning.
\end{IEEEbiography}
\begin{IEEEbiography}[{\includegraphics[width=1in,height=1.25in,clip,keepaspectratio]{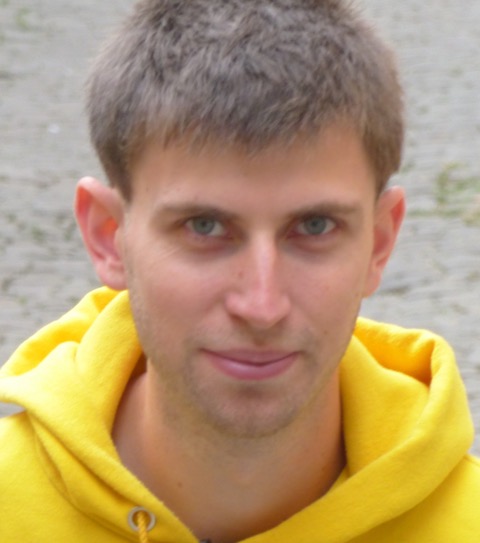}}]{Torbjørn Cunis} (M'2022) received the M.Sc. degree in automation
engineering from the RWTH Aachen University in 2016
and his doctoral degree in systems and control from
ISAE-Supaéro, Toulouse, in 2019. Since 2021, he has been a lecturer at the Institute of Flight Mechanics and Controls of the University of
Stuttgart.  He was a research fellow of the University of Michigan from 2019 to 2021. His research focuses on systems theory of nonconvex optimization algorithms, optimization-based analysis and control methods for nonlinear dynamical systems, and nonlinear sum-of-squares optimization. Dr. Cunis is an adjunct of the University of Michigan Aerospace Department.
\end{IEEEbiography}




\end{document}